\documentclass[aop,seceqn,dvips]{arximspdf}

\doi{10.1214/009117906000001178}
\volume{35}
\issue{6}
\pubyear{2007}
\firstpage{2213}
\lastpage{2262}

\makeatletter
\newtheorem{theo}{Theorem}[section]
\newtheorem{lemma}[theo]{Lemma}
\newtheorem{prop}[theo]{Proposition}
\newproclaim{rem}[theo]{Remark}
\newtheorem{cor}[theo]{Corollary}
\newproclaim{defi}[theo]{Definition}
\newproclaim{exemple}[theo]{Example}

\newtheorem{ta}{Technical Assumption}
\makeatother

\begin{document}
\begin{frontmatter}

\title{Some parabolic PDE\textup{s} whose drift is an irregular~random~noise~in~space}
\runtitle{Some parabolic PDE\textup{s}}

\begin{aug}
\author[A]{\fnms{Francesco}
\snm{Russo}\corref{}\ead[label=e1]{russo@math.univ-paris13.fr}} and
\author[B]{\fnms{Gerald} \snm{Trutnau}\ead[label=e2]{trutnau@math.uni-bielefeld.de}}
\runauthor{F. Russo and G. Trutnau}
\affiliation{Universit\'{e} Paris \textit{13} and Universit\"at Bielefeld}
\address[A]{Institut Galil\'{e}e, Math\'ematiques\\
Universit\'{e} Paris 13\\
99, Avenue J.B.~Cl\'{e}ment\\
F-93430 Villetaneuse\\
France\\
\printead{e1}} 
\address[B]{Universit\"at Bielefeld\\
Fakult\"at f\"ur Mathematik\\
Postfach 10 01 31 \\
D-33501 Bielefeld\\
Germany\\
\printead{e2}}
\end{aug}

\received{\smonth{3} \syear{2006}}
\revised{\smonth{8} \syear{2006}}

%
\begin{abstract}
A new class of random partial differential equations of parabolic type is
considered, where the stochastic term consists of an irregular noisy drift,
not necessarily Gaussian, for which a suitable interpretation is provided.
After freezing a realization of the drift (stochastic process),
we study existence and uniqueness (in some appropriate sense)
of the associated parabolic equation and a probabilistic interpretation
is investigated.

\end{abstract}

%

\end{frontmatter}

\section{Introduction}\label{s1}

This paper focuses on a random partial differential equation consisting
of a parabolic PDE with irregular noise in the drift. Formulation,
existence (with uniqueness in a certain sense) and double probabilistic
representation are discussed. The equation itself is motivated by
\textit{random irregular media models}.

Let $T > 0$, $ \sigma\dvtx {\mathbb R}\rightarrow{\mathbb R}$ be a
continuous function and $\dot\eta(x)$ a generalized random field
playing the role of a noise. Let $u^0 \dvtx {\mathbb
R}\rightarrow{\mathbb R}$, $\lambda\dvtx  [0,T] \times{\mathbb
R}\rightarrow{\mathbb R}$ be continuous. Consider the problem
\begin{eqnarray}\label{eI-0}
-\partial_t v(t,x) + \frac{\sigma^2(x)}{2}\,
\partial^2_{\mathit{xx}} v(t,x) +
\dot\eta(x)\, \partial_{x} v(t,x) &=& \lambda(T-t,x),
\nonumber
\\[-8pt]
\\[-8pt]
\nonumber v(0,x) &=& u^0 (x),
\end{eqnarray}
where $\dot\eta$ is the derivative in the sense of distributions of a
continuous process. Among examples of possible $\eta$, we have in mind
not only different possibilities of continuous processes as classical
Wiener process and (multi) fractional Brownian motion, but also
non-Gaussian processes. The derivative in the sense of distributions
$\dot\eta(x)$ will be the associated noise. (\ref{eI-0}) is a new type
of SPDE, not yet studied in any real depth even when $\eta$ is a classical
Brownian motion. For the situation where $\dot\eta(x)$ is replaced by a
space-time white noise $\dot\eta(t,x)$, some relevant work was done by
Nualart and Viens (see, e.g., \cite{nv}). In this article, time dependence
is useful for the corresponding stochastic integration.

Equation~(\ref{eI-0}) is equivalent to the following \textit{dual}
problem:
\begin{eqnarray}\label{eI-0dual}
\partial_t u(t,x) + \frac{\sigma^2(x)}{2}
\,\partial^2_{\mathit{xx}} u(t,x) + \dot\eta(x) \,\partial_{x} u(t,x)
&=& \lambda(t,x), \nonumber
\\[-8pt]
\\[-8pt]
\nonumber u(T,x) &=& u^0 (x).
\end{eqnarray}
Formally speaking, setting $u(t,x) = v(T-t,x)$, $v$ solves (\ref{eI-0})
if and only if $u$ solves (\ref{eI-0dual}). This is rigorously
confirmed in Section~\ref{s9} so that at this stage, the choice of
whether to work with equation (\ref{eI-0}) or (\ref{eI-0dual}) is
arbitrary. We have decided to concentrate on equation (\ref{eI-0dual})
because it corresponds to the standard form for probabilistic
representation.

The idea of this paper is to first freeze the realization $\omega$, to
set $b(x) = \eta(x) (\omega)$ and then to consider the deterministic
Cauchy problem associated with~(\ref{eI-0dual}),
\begin{eqnarray}\label{eI-0det}
\partial_t u(t,x) + \frac{\sigma^2(x)}{2}
\,\partial^2_{\mathit{xx}} u(t,x) + b' (x)\, \partial_{x} u(t,x) &=&
\lambda(t,x), \nonumber
\\[-8pt]
\\[-8pt]
\nonumber u(T,x) &=& u^0 (x),
\end{eqnarray}
where $b'$ is the derivative of the continuous function $b$.

Since the product of a distribution and a continuous function is not
defined in the theory of Schwarz distributions, we must develop some
substitution tools. Ideally, we would like to represent the parabolic
PDE probabilistically through a diffusion which is the solution of the
stochastic differential equation (SDE)
\begin{equation} \label{eI-2}
dX_t = \sigma(X_t)\,dW_t + b'(X_t)\,dt
\end{equation}
with generalized drift. We will give a meaning to (\ref{eI-2}) at three
different levels:
\begin{itemize}
\item the level of a martingale problem;

\item the level of a stochastic differential equation in the sense of
probability laws;

\item the level of a stochastic differential equation in the strong
sense.
\end{itemize}
For each of these levels, we shall provide conditions for equation
(\ref{eI-2}), with given initial data, to be well posed. Later, the
notion of a $C^0_b$-solution to the generalized parabolic PDE
(\ref{eI-0det}) will be defined; related to this, existence, uniqueness
and probabilistic representation will be shown.

When $\eta$ is a strong finite cubic variation process and $\sigma= 1$,
the solutions to~(\ref{eI-0det}) obtained for $b = \eta(\omega)$
provide solutions to the SPDE (\ref{eI-0}). This is shown in the last
part of the paper. A typical example of a strong zero cubic variation
process is the fractional Brownian motion with Hurst index $H
\ge\frac{1}{3}$. Equation (\ref{eI-0det}) will be understood in some
\textit{weak distributional} sense that we can formally reconstruct as
follows. We freeze $b = \eta(\omega)$ as a realization and formally
integrate equation (\ref{eI-0}) from $0$ to $t$ in time against a
smooth test function $\alpha$ with compact support in space. The result
is
\begin{eqnarray} \label{eI-1}
\nonumber \hspace*{10mm}  && - \int_{\mathbb R}\, dx\, \alpha(x) u(t,x)
+ \int_{\mathbb R} dx\, \alpha(x) u^0(x) - \int_0^t \,ds\, \tfrac{1}{2}
\int_{\mathbb R}dx\, \alpha' (x)\, \partial _x u(s,x)
\\
&&\quad{} + \int_0^t ds \int_{\mathbb R}b (dx) \alpha(x)\, \partial_x
u(s,x)
\\
\nonumber &&\qquad = \int_0^t ds \int_{\mathbb R}dx\, \alpha(x)
\lambda(T-s,x).
\end{eqnarray}
The integral $\int_{\mathbb R}\alpha(x)\, \partial_x u (s,x) b(dx)$
needs interpretation since $b$ is not generally of bounded variation
and it involves the product of the distribution $b'$ and the function
$\partial_x u (s,\cdot)$; in general, this function is, unfortunately,
only continuous. As expected this operation is deterministically
undefined, unless one uses a generalized functions theory. However,
since $b$ is a frozen realization of a stochastic process $\eta$, we
can hope to justify the integral in a stochastic sense. Note that it
cannot be of It\^o type, even if $\eta$ were a semimartingale, since
$\partial_x u(s,\cdot)$ is not necessarily adapted to some
corresponding filtration. We will, in fact, interpret the stochastic
integral element $b(dx)$ or $\eta(dx)$ as a symmetric (Stratonovich)
integral $d^0\eta$ of regularization type; see Section~\ref{s3}.

\begin{defi} \label{deI-1rig}
A continuous random field $(v(t,x), t \in[0,T], x \in{\mathbb R}) $,
a.s. in $ C^{0,1} (]0,T[ \times{\mathbb R})$, is said to be a
(\textit{weak}) \textit{solution} to the SPDE (\ref{eI-0}) if
\begin{eqnarray} \label{eI-1rig}
\nonumber \hspace*{10mm} && - \int_{\mathbb R} dx\, \alpha(x) v(t,x) +
\int_{\mathbb R} dx\, \alpha(x) v^0(x) - \int_0^t ds\, \tfrac{1}{2}
\int_{\mathbb R}dx\, \alpha' (x)\, \partial_x v(s,x)
\\
&&\quad{} + \int_{\mathbb R} d^\circ\eta(x) \alpha(x)  \biggl( \int_0^t
ds\,
\partial_x v(s,x)  \biggr)
\\
\nonumber &&\qquad = \int_0^t ds \int_{\mathbb R}dx\, \alpha(x)
\lambda(T-s,x)
\end{eqnarray}
for every smooth function with compact support $\alpha$.
\end{defi}

If we integrate equation (\ref{eI-0dual}) from $t$ to $T$ in time
against a smooth test function~$\alpha$ with compact support in space,
we are naturally led to the following.

\begin{defi} \label{deI-1dual}
A continuous random field $(u(t,x), t \in[0,T], x \in{\mathbb R}) $,
a.s. in $ C^{0,1} (]0,T[ \times{\mathbb R})$, is said to be a
(\textit{weak}) \textit{solution} to the SPDE (\ref{eI-0dual}) if
\begin{eqnarray} \label{eI-1dual}
\nonumber \hspace*{10mm} && - \int_{\mathbb R} dx\, \alpha(x) u(t,x) +
\int_{\mathbb R} dx\, \alpha(x) u^0(x) - \int_t^T ds\, \tfrac{1}{2}
\int_{\mathbb R}dx\, \alpha' (x)\, \partial_x u(s,x)
\\
&&\quad{} + \int_{\mathbb R}d^\circ\eta(x) \alpha(x)  \biggl( \int_t^T
ds\, \partial_x u(s,x)  \biggr)
\\
\nonumber &&\qquad = \int_t^T ds \int_{\mathbb R}dx\, \alpha(x)
\lambda(s,x)
\end{eqnarray}
for every smooth function with compact support $\alpha$.
\end{defi}

We will show that the \textit{probabilistic} solutions that we
construct through stochastic equation (\ref{eI-2}) will, in fact, solve
(\ref{eI-1}).

Diffusions in the generalized sense were studied by several authors
beginning with (at least to our knowledge) \cite{po1}. Later, many
authors considered special cases of stochastic differential equations
with generalized coefficients.  It is difficult to quote them all. In
particular, we refer to the case when $b$ is a measure
\cite{esd,made,o}. In all of these cases, solutions were
semimartingales. More recently, \cite{ew} considered special cases of
nonsemimartingales solving stochastic differential equations with
generalized drift; those cases include examples coming from Bessel
processes.

\cite{frw1} and \cite{frw2} treated well-posedness of the martingale
problem, It\^o's formula under weak conditions, semimartingale
characterization and the Lyons--Zheng decomposition. The only
assumption was the strict positivity of $\sigma$ and the existence of
the function $\Sigma(x) = 2 \int_0^x \frac{b'}{\sigma ^2}\, dy$ with
appropriate regularizations. Bass and Chen~\cite{BC} were also
interested in (\ref{eI-2}) and provided a well-stated framework
when $\sigma$ is $\frac{1}{2}$-H\"older continuous and $b$ is
$\gamma$-H\"older continuous, $ \gamma> \frac{1}{2}$.

Beside the martingale problem, in the present paper, we shall emphasize
the formulation of (\ref{eI-2}) as a stochastic differential equation
which can be solved by introducing more assumptions on the
coefficients. Several examples are provided for the case of weak and
strong solutions of (\ref{eI-2}).

The paper is organized as follows. Section~\ref{s2} is devoted to basic
preliminaries, including definitions and properties related to Young
integrals. Section~\ref{s3} is devoted to some useful remainder in
stochastic calculus via regularization. In Section~\ref{s4}, we
introduce the formal \textit{elliptic} operator $L$ and recall the
concept of a $C^1$-generalized solution of $Lf = \dot\ell$ for
continuous real functions $\dot\ell$. We further introduce a
fundamental hypothesis on $L$ for the sequel, called Technical
Assumption~\textup{\ref{A0}}, and we illustrate several examples where
it is verified. In Section~\ref{s5}, we discuss different notions of
martingale problems. Section~\ref{s6} provides notions of solutions to
\textit{stochastic differential equations with distributional drift}
and their connections with martingale problems. The notion of solution
is coupled with a property of \textit{extended local time regularity}.
This concept of solution is new, even when the drift is an ordinary
function. Section~\ref{s7} presents the notion of a $C^0_b$-solution
for a parabolic equation ${\mathcal  L}u = \lambda$, where $\lambda$ is
bounded and continuous with ${\mathcal  L}=
\partial_t + L$. We also provide existence, uniqueness and
probabilistic representations of $C^0_b$-solutions to ${\mathcal  L}u =
\lambda$. Section~\ref{s8} discusses \textit{mild} solutions to the
previous parabolic PDE and useful integrability properties for its
solutions. In Section~\ref{s9}, we finally show that the
$C^0_b$-solutions provide, in fact, true \textit{weak} solutions to the
SPDE (\ref{eI-0}) if $\sigma= 1$.

\section{Preliminaries}\label{s2}

In this paper, $T$ will be a fixed horizon time, unless otherwise
specified. A function $f$ defined on $[0,T]$ (resp., ${\mathbb R}_+$)
will be extended, without mention, by setting $f(t) = f(0)$ for $t
\le0$ and $f(T)$ for $t \ge T$ [resp., $f(0)$ for $t \le0$].

$C^0({\mathbb R}) $ will indicate the set of continuous functions
defined on ${\mathbb R}, $ $C^p({\mathbb R}),$ the space of real
functions with differentiability class $C^p$. We denote by
$C^0_0({\mathbb R})$ [resp., $C^1_0({\mathbb R})$] the space of
continuous (continuous differentiable) functions vanishing at zero.
When there is no confusion, we will also simply use the symbols $C^0,
C^p, C^0_0, C^1_0 $. We denote by $ C^0_b ([0,T] \times{\mathbb R})$
the space of real continuous bounded functions defined on $[0,T] \times
{\mathbb R}$. $C^0_b ({\mathbb R})$, or simply $C^0_b $, indicates the
space of continuous bounded functions defined on ${\mathbb R}$.

The vector spaces $C^0 ({\mathbb R})$ and $C^p ({\mathbb R})$ are
topological Fr\'{e}chet spaces, or \mbox{F-}spaces, according to the
terminology of \cite{ds}, Chapter~1.2. They are equipped with the
following natural topology. A sequence $f_n$ belonging to $C^0({\mathbb
R})$ [resp., $C^p ({\mathbb R})$] is said to converge to $f$ in the
$C^0({\mathbb R})$ [resp., $C^p({\mathbb R})$] sense if $f_n$ (resp.,
$f_n$ and all derivatives up to order $p$) converges (resp.,
converge) to $f$ (resp., to $f$ and all its derivatives) uniformly on
each compact of ${\mathbb R}$.

We will consider functions $ u \dvtx [0,T] \times{\mathbb R}
\rightarrow{\mathbb R}$ which are bounded and continuous. A sequence
$(u_n)$ in $C^0_b([0,T] \times{\mathbb R})$ will be said to converge
\textit{in a bounded way} to $u$ if:

\begin{itemize}
\item $\lim_{n \to\infty} u_n(t,x) = u(t,x), \quad\forall(t,x) \in[0,T]
\times{\mathbb R}$;

\item there exists a constant $ c>0$, independent of the sequence, such
that
\begin{equation}\label{e1-1}
\sup_{t \le T, x \in{\mathbb R}} | u_n(t,x) |\le c\qquad \forall n
\in{\mathbb N}.
\end{equation}
\end{itemize}
If the sequence $(u_n)$ does not depend on $t$, we similarly define the
convergence of $(u_n) \in C^0_b({\mathbb R})$ to $u \in C^0_b({\mathbb
R})$ in a bounded way.

Given two functions $ u_1, u_2\dvtx  [0,T] \times{\mathbb R}
\rightarrow {\mathbb R}$, the composition notation $ u_1 \circ u_2$
means $ (u_1 \circ u_2)(t,x) = u_1(t,u_2(t,x))$.

For positive integers $m,k$, $C^{m,k} $ will indicate functions in the
corresponding differentiability class. For instance, $ C^{1,2} ([0,T[
\times{\mathbb R})$ will be the space of $ (t,x) \mapsto u(t,x)$
functions which are $C^1$ on $[0,T[ \times{\mathbb R}$ (i.e., once
continuously differentiable) and such that $ \partial_{\mathit{xx}}^2
u$ exists and is continuous.

$C^{m,k}_b $ will indicate the set of functions $C^{m,k} $ such that
the partial derivatives of all orders are bounded.

If $I$ is a real compact interval and $\gamma\in\,]0,1[$, we denote by $
C^\gamma(I)$ the vector space of real functions defined on $I$ which
are H\"older with parameter $\gamma$. We denote by $ C^\gamma({\mathbb
R})$, or simply $C^\gamma$, the space of locally H\"older functions,
that is, H\"older on each real compact interval.

Suppose $I = [\tau, T]$, $\tau, T$ being two real numbers such that
$\tau< T$. Here, $T$ does not necessarily need to be positive. Recall
that $f \dvtx  I \mapsto\mathbb{R}$ belongs to $C^\gamma(I)$ if
\[
N_\gamma(f):=\sup_{\tau\leq s,t \leq
T}\frac{|f(t)-f(s)|}{|t-s|^\gamma}<\infty.
\]
Clearly, $f \mapsto| f(\tau) |+ N_\gamma(f) $ defines a norm on
$C^\gamma(I)$ which makes it a Banach space. $ C^\gamma({\mathbb R})$
is an F-space if equipped with the topology of convergence related to $
C^\gamma(I)$ for each compact interval $I$. A sequence $(f_n)$ in
$C^\gamma({\mathbb R})$ converges to $f$ if it converges according to
$C^\gamma(I)$ for every compact interval $I$.

We will also provide some reminders about the so-called \textit{Young
integrals} (see~\cite{y}) but will remain, however, in a simplified
framework, as in \cite{fepr} or \cite{RVSem}. We recall the essential
inequality, stated, for instance, in \cite{fepr}:

Let $ \gamma, \beta> 0$ be such that $\gamma+ \beta> 1$. If $f, g \in
C^1(I)$, then
\begin{equation} \label{Fyoung}
\bigg|\int_a^b \bigl(f(x) - f(a)\bigr)\,dg(x)  \bigg| \leq
C_\rho(b-a)^{1+\rho } N_\gamma(f) N_\beta(g)
\end{equation}
for any $[a,b]\subset I$ and $\rho\in\,]0,\gamma+\beta-1[$, where
$C_\rho$ is a constant not depending on $f,g$. The bilinear map sending
$(f,g)$ to $\int_0^{\bolds{\cdot}} f\, dg$ can be continuously extended
to $C^\gamma(I) \times C^\beta(I) $ with values in $C^0(I) $. By
definition, that object will be called the \textit{Young integral} of
$f$ with respect to $g$ on $I$. We also denote it $
\int_\tau^{\bolds{\cdot}} f\, d^{(y)} g $.

By additivity, we set, for $a, b \in[\tau,T]$,
\[
\int_ a^b f\, d^{(y)} g = \int_ \tau^b f\, d^{(y)} g - \int_ \tau^a f\,
d^{(y)} g.
\]

Moreover, the bilinear map defined on $C^1({\mathbb R}) \times
C^1({\mathbb R})$ by $(f, g) \rightarrow\int_0^{\bolds{\cdot}} f\, dg$
extends continuously to $C^\gamma({\mathbb R}) \times C^\beta ({\mathbb
R})$ onto $C^0({\mathbb R})$. Again, that object, defined on the whole
real line, will be called \textit{Young integral} of $f$ with respect
to $g$ and will again be denoted by $ \int_0^{\bolds{\cdot}} f\,
d^{(y)} g $.

\begin{rem} \label{REXVQ7}
Inequality (\ref{Fyoung}) remains true for $f \in C^\gamma(I), g \in
C^\beta(I)$. In particular, $t \mapsto\int_\tau^t f\, d^{(y)} g$
belongs to $ C^\beta(I)$. In fact,
\[
\bigg|\int_a^b f\, dg  \bigg| \le \bigg|\int_a^b \bigl(f -
f(a)\bigr)\,dg  \bigg| +  \big| f(a)\bigl( g(b) - g(a)\bigr) \big|.
\]
\end{rem}

Through the extension of the bilinear operator sending $(f,g)$ to $\int
_0^{\bolds{\cdot}} f\, dg$, it is possible to get the following chain
rule for Young integrals.

\begin{prop} \label{PCRuleY}
Let $f, g, F\dvtx  I \rightarrow{\mathbb R}$, $I = [\tau, T]$. We
suppose that $g \in C^\beta(I)$, $f \in C^\gamma(I)$, $F \in
C^\delta(I)$ with $\gamma+ \beta> 1$, $\delta+ \beta> 1.$ We define
$G(t) = \int_\tau^t f\, d^{(y)} g$. Then
\[
\int_\tau^t F\, d^{(y)} G = \int_\tau^t F f\, d^{(y)} g.
\]
\end{prop}

\begin{pf}
If $g \in C^1(I)$, then the result is obvious. We remark that $G \in
C^\gamma(I)$. Repeatedly using inequality (\ref{Fyoung}), one can
show that the two linear maps $g \mapsto\int_\tau^{\bolds{\cdot}} F\,
d^{(y)} G$ and $g \mapsto\int_\tau^{\bolds{\cdot}} F f\, d^{(y)} g$ are
continuous from $C^\delta(I) $ to $C^0(I)$. This concludes the proof of
the proposition.
\end{pf}

By a \textit{mollifier}, we mean a function $\Phi\in{\mathcal
{S}}(\mathbb{R})$ (i.e.,  a $C^{\infty}$-function such that itself and
all its derivatives decrease to zero faster than any power of
$|x|^{-1}$ as $|x|\to\infty$) with $\int\Phi(x)\,dx = 1$. We set
$\Phi_n (x) := n \Phi(nx)$.

The result below shows that mollifications of a H\"older function $f$
converge to $f$ with respect
to the H\"older topology.

\begin{prop}\label{PYoungreg}
Let $\Phi$ be a mollifier and let $f \in C^{\gamma'} (I)$. We write
$f_n = \Phi_n \ast f$. Then $f_n \rightarrow f$ in the $C^\gamma(I)$
topology for any $0 < \gamma< \gamma'$.
\end{prop}

\begin{pf}
We need to show that $N_\gamma(f - f_n)$ converges to zero. We set
$\Delta_n (t) = (f - f_n) (t)$. Let $a, b \in I$. We will establish
that
\begin{equation} \label{NYoung}
|\Delta_n (b) - \Delta_n (a) |\le\operatorname{const} | b - a |^\gamma
\biggl( \frac{1}{n}  \biggr)^{\gamma' - \gamma}.
\end{equation}
Without loss of generality, we can suppose that $a < b $. We
distinguish between two cases.

\textit{Case} $ a < a + \frac{1}{n} < b$.

\qquad We have
\begin{eqnarray*}
|\Delta_n (b) - \Delta_n (a) |& \le&
 \bigg|\int\biggl( f\biggl(b - \frac{y}{n}\biggr) - f(b)\biggr) \Phi(y)\,dy  \bigg|
\\
&&{} +  \bigg|\int\biggl( f\biggl(a - \frac{y}{n}\biggr) - f(a)\biggr)
\Phi(y)\,dy \bigg|
\\
& \le& 2 \int \bigg|\frac{y}{n} \bigg|^{\gamma'} |
\Phi(y) |\,dy
\\
&\le& 2 \int |\Phi(y) || y |^{\gamma'} \,dy (b-a)^\gamma
\biggl(\frac{1}{n} \biggr)^{\gamma' - \gamma}.
\end{eqnarray*}

\textit{Case} $ a < b \le a + \frac{1}{n} $.

\qquad In this case, we have
\begin{eqnarray*}
&& |\Delta_n (b) - \Delta_n (a) |
\\
&&\qquad  \le \int| f(b) - f(a) ||\Phi(y) |\,dy + \int\bigg| f\biggl(b
+ \frac{y}{n} \biggr) - f\biggl(a + \frac{y}{n}\biggr) \bigg||\Phi
(y) |\,dy
\\
&&\qquad \le 2 ( b - a) ^{\gamma'}\int|\Phi(y) |\,dy  \le2 \int|\Phi(y)
|\,dy (b-a)^\gamma \biggl(\frac {1}{n} \biggr)^{\gamma' - \gamma}.
\end{eqnarray*}

Therefore, (\ref{NYoung}) is verified with $\operatorname{const} = 2
\int|\Phi (y) |(1+| y |^{\gamma'} )\,dy$. This implies that
\[
N_\gamma(f - f_n) \le\operatorname{const}  \biggl(\frac{1}{n} \biggr)^
{\gamma' - \gamma},
\]
which allows us to conclude.
\end{pf}

For convenience, we introduce the topological vector space defined by
\[
D^\gamma= \bigcup_{\gamma' > \gamma} C^{\gamma'}({\mathbb R}).
\]
It is also a \textit{vector algebra}, that is, $D^\gamma$ is a vector
space and an algebra with respect to the \textit{sum and product of
functions}.

The next corollary is a consequence of the definition of the Young
integral and Remark~\ref{REXVQ7}.

\begin{cor}\label{CYoung}
Let $ f \in D^\gamma$, $g \in D^\beta$ with $\gamma+ \beta\ge1$. Then
$t \mapsto\int_0^t f\, d^{(y)} g $ is well defined and belongs to
$D^\beta$.
\end{cor}

$D^\gamma$ is not a metric space, but an inductive limit of the
F-spaces $C^\gamma$; the weak version of the Banach--Steinhaus theorem
for F-spaces can be adapted.

In fact, a direct consequence of the Banach--Steinhaus theorem of
\cite{ds}, Section~2.1, is the following.

\begin{theo} \label{TBS}
Let $E = \bigcup_n E_n$ be an inductive limit of F-spaces $E_n$ and $F$
another F-space. Let $(T_n)$ be a sequence of continuous linear
operators \mbox{$T_n \dvtx  E \rightarrow F$}. Suppose that $ T f :=
\lim_{n \rightarrow\infty} T_n f$ exists for any $f \in E$. Then
\mbox{$T \dvtx  E \rightarrow F$} is again a continuous (linear)
operator.
\end{theo}

\section{Previous results in stochastic calculus via
regularization}\label{s3}

We recall here a few notions related to stochastic calculus via
regularization, a theory which began with~\cite{rv1}. We refer to a
recent survey paper \cite{RVSem}.

The stochastic processes considered may be defined on $[0,T], {\mathbb
R}_+$ or ${\mathbb R}$. Let $X=(X_t, t \in{\mathbb R})$ be a continuous
process and $Y=(Y_t, t \in{\mathbb R})$ be a process with paths in
$L^1_{\mathrm{loc}}$. For the paths of process $Y$ with parameter on
$[0,T]$ (resp., ${\mathbb R}_+$), we apply the same convention as was
applied at the beginning of previous section for functions. So we
extend them without further mention, setting $Y_0$ for $t \le0$ and $Y_T$
for $t \ge T$ (resp., $Y_0$ for $t \le0$). ${\mathcal  C}$ will denote
the vector algebra of continuous processes. It is an F-space if
equipped with the topology of u.c.p. (uniform convergence in
probability) convergence.

In the sequel, we recall the most useful rules of calculus; see, for
instance, \cite{RVSem} or \cite{rv4}.

The forward symmetric integrals and the covariation process are defined
by the following limits in the u.c.p. sense, whenever they exist:
\begin{eqnarray}
\label{E1.1} \int_0^t Y\,d^{-} X & := & \lim_{\varepsilon\to0+}
\int_0^t Y_s \frac{X_{s+\varepsilon}-X_s}{\varepsilon}\,   ds,
\\
\label{E1.1S} \int_0^t Y_s \, d^{\circ} X_s & := & \lim_{\varepsilon
\to0+} \int_0^t Y_s \frac{X_{s+\varepsilon}-X_{s -\varepsilon} }{2
\varepsilon} \,  ds,
\\
\label{E1.4}  [X,Y ]_t & := & \lim_{\varepsilon\to0+}
C^{\varepsilon}(X,Y)_t,
\end{eqnarray}
where
\[
C^{\varepsilon}(X,Y)_t := \frac{1}{\varepsilon} \int_0^t
(X_{s+\varepsilon} - X_s) (Y_{s+\varepsilon}-Y_s)\,ds.
\]
All stochastic integrals and covariation processes will of course be
elements of~$\mathcal {C}$. If $[X,Y]$, $[X,X]$ and $[Y,Y]$ exist, we
say that $(X,Y)$ \textit{has all of its mutual covariations}.

\begin{rem} \label{R1.0} If $X$ is (locally) of bounded variation, we
have:
\begin{itemize}
\item $ \int_0^t X \, d^{-} Y = \int_0^t X_s \, d^{\circ} Y_s =
\int_0^t X_s \, dY_s $, where the third integral is meant in the
Lebesgue--Stieltjes sense;

\item $ [X,Y] \equiv0$.
\end{itemize}
\end{rem}

\begin{rem} \label{R1.1}
(a) $ \int_0^t Y_s \, d^{\circ} X_s = \int_0^t Y_s \, d^{-} X_s +
\frac{1}{2} [X,Y]$ provided that two of the three integrals or
covariations exist.
{\smallskipamount=0pt
\begin{longlist}[(a)]
\item[(b)] $X_t Y_t = X_0 Y_0 + \int_0^t Y_s \, d^{-} X_s + \int_0^t
X_s\, d^{-} Y_s + [X,Y]_t $ provided that two of the three integrals or
covariations exist.

\item[(c)] $X_t Y_t = X_0 Y_0 + \int_0^t Y\,  d^{\circ} X + \int_0^t
X_s \,d^{\circ} Y_s$ provided that one of the two integrals exists.
\end{longlist}}
\end{rem}

\begin{rem} \label{R1.2}
(a) If $[X,X]$ exists, then it is always an increasing process and $X$
is called a \textit{finite quadratic variation process}. If $[X,X] =
0$, then $X$ is said to be a \textit{zero quadratic variation process}.
{\smallskipamount=0pt
\begin{longlist}[(a)]
\item[(b)] Let $X$, $Y$ be continuous processes such that $(X,Y)$ has
all of its mutual covariations. Then $[X,Y]$ has locally bounded
variation. If $f,g\in C^1$, then
\[
[f(X),g(Y)]_t = \int_0^t f'(X)g'(Y)\,d[X,Y].
\]

\item[(c)] If $A$ is a zero quadratic variation process and $X$ is a
finite quadratic variation process, then $[X,A] \equiv0$.

\item[(d)] A bounded variation process is a zero quadratic variation
process.

\item[(e)] (\textit{Classical It\^o formula}.) If $f\in C^2$, then
$\int_0^{\bolds{\cdot}} f'(X)\,d^-X$ exists and is equal to
\[
f(X)- f(X_0) - \tfrac{1}{2} \int_0^{\bolds{\cdot}} f''(X)\,d[X,X].
\]

\item[(f)] If $g \in C^1$ and $f \in C^2$, then the forward integral $
\int_0^{\bolds{\cdot}} g(X)\,d^- f(X)$ is well defined.
\end{longlist}}
\end{rem}

In this paper, all filtrations are supposed to fulfill the usual
conditions. If $\mathbb{F} = ({\mathcal {F}}_t)_{t\in[0,T]}$ is a
filtration, $X$ an $\mathbb{F}$-semimartingale and $Y$ is an
$\mathbb{F}$-adapted cadlag process, then $\int_0^{\bolds{\cdot}} Y\,
d^{-}X$ is the usual It\^{o} integral. If $Y$ is an
$\mathbb{F}$-semimartingale, then $\int_0^{\bolds{\cdot}} Y \,  d^{\circ}
X$ is the classical Fisk--Stratonovich integral and $[X,Y]$ is the
usual covariation process $\langle X,Y\rangle$.

We now introduce the notion of Dirichlet process, which was essentially
introduced by F\"ollmer \cite{fodir} and has been considered by many
authors; see, for instance, \cite{ber,rvw} for classical properties.

In the present section, $(W_t)$ will denote a classical $({\mathcal
F}_t)$-Brownian motion.

\begin{defi} \label{D81}
An $({\mathcal  F}_t)$-adapted (continuous) process is said to be a
\mbox{$({\mathcal F}_t)$-}\textit{Dirichlet process} if it is the sum
of an \mbox{$({\mathcal F}_t)$-}local martingale $M$ and a zero quadratic
variation process $A$. For simplicity, we will suppose that $A_0 = 0$
a.s.
\end{defi}

\begin{rem}\label{R81}
(i) Process $(A_t)$ in the previous decomposition is an\break \mbox{$({\mathcal
F}_t)$-}adapted process.
{\smallskipamount=0pt
\begin{longlist}[(iii)]
\item[(ii)]An $({\mathcal  F}_t)$-semimartingale is an $({\mathcal
F}_t)$-Dirichlet process.

\item[(iii)] The decomposition $M + A$ is unique.

\item[(iv)] Let $f\dvtx  {\mathbb R}\rightarrow{\mathbb R}$ be of
class~$C^1$ and let $X$ be an $({\mathcal  F}_t)$-Dirichlet process.
Then $f(X)$ is again an $({\mathcal F}_t)$-Dirichlet process with local
martingale part $M^f_t = f(X_0) + \int_0^t f'(X)\,dM$.
\end{longlist}}
\end{rem}

The class of semimartingales with respect to a given filtration is
known to be stable with respect to $C^2$ transformations.
Remark~\ref{R1.2}(b) says that finite quadratic variation processes are
stable through $C^1$ transformations. The last point of the previous
remark states that $C^1$ stability also holds for Dirichlet processes.

Young integrals introduced in Section~\ref{s2} can be connected with
the forward and symmetric integrals via the regularization appearing
before Remark~\ref{R1.0}. The next proposition was proven
in~\cite{RVSem}.

%
\begin{prop} \label{PEXVQ7}
Let $ X, Y$ be processes whose paths are respectively in $ C^\gamma$
and $ C^\beta$, with $\gamma
>0$, $\beta>0$ and $\gamma+ \beta>1$.

For any symbol $\star\in\{-,\circ\}$, the integral $
\int_0^{\bolds{\cdot}} Y\,d^\star X$ coincides with the Young integral
$ \int_0^{\bolds{\cdot}} Y\,d^{(y)} X$.
\end{prop}
%


\begin{rem} \label{REXVQ7a}
Suppose that $X$ and $Y$ satisfy  the conditions of
Proposition~\ref{PEXVQ7}. Then Remark~\ref{R1.1}(a) implies that
$[X,Y]=0$.
\end{rem}

We need an extension of stochastic calculus via regularization in the
direction of higher $n$-variation. The properties concerning variation
higher than $2$ can be found, for instance, in \cite{er}.

We set
\[
[X,X,X]^\varepsilon_t =\frac{ 1}{\varepsilon} \int_0^t
(X_{s+\varepsilon}-X_s)^3\, ds.
\]
We also define
\[
 \|[X,X,X]^\varepsilon \|_t = \frac{1}{\varepsilon} \int_0^t
| X_{s+\varepsilon}-X_s|^3 \, ds.
\]
If the limit in probability of $[X,X,X]^\varepsilon_t$ when
$\varepsilon\rightarrow0$ exists for any $t$, we denote it by
$[X,X,X]_t$. If the limiting process $[X,X,X]$ has a continuous
version, we say that $X$ is a \textit{finite cubic variation process}.

If, moreover, there is a positive sequence $(\varepsilon_n)_{n\in
{\mathbb N}}$ converging to zero such that
\begin{equation}
\sup_{\varepsilon_n}  \|[X,X,X]^{\varepsilon_n}  \|_T < +\infty,
\end{equation}
then we say that $X$ is a (\textit{strong}) \textit{finite cubic variation process}.
If $X$ is a (strong) finite cubic variation process such that $[X,X,X]
= 0$, then $X$ will be said to be a (strong) zero \textit{finite cubic
variation process}.

For instance, if $X = B^H$, a fractional Brownian motion with Hurst
index $H$, then $X$ is a finite quadratic variation process if and only if
$H \ge \frac{1}{2}$; see \cite{rv4}. It is a strong zero cubic
variation process if and only if $H \ge\frac {1}{3}$; see \cite{er}. On
the other hand, $B^H$~is a zero cubic variation process if and only if
$H > \frac{1}{6}$; see \cite{GRV}.

It is clear that a finite quadratic variation process is a strong zero
cubic variation process. On the other hand, processes whose paths are
H\" older continuous with parameter greater than $\frac{1}{3}$ are
strong zero cubic variation processes.

As for finite quadratic variation and Dirichlet processes, the
$C^1$-stability also holds for finite cubic variation processes. The
next proposition is a particular case of a result contained in
\cite{er}.

\begin{prop} \label{Pstacub} Let $X$ be a strong finite cubic
variation process, $V$ a locally bounded variation process and $f\dvtx
{\mathbb R}\times{\mathbb R}\rightarrow{\mathbb R}$ of class $C^{1}$.
Then $Z = f(V,X)$ is again a strong finite cubic variation process and
\[
[Z,Z,Z]_t = \int_0^t \partial_x f(V_s,X_s)^3\, d[X,X,X]_s.
\]
\end{prop}

Moreover, an It\^o chain rule property holds, as follows.

\begin{prop} \label{Crule} Let $X$ be a strong finite cubic variation process,
$ V$ a bounded variation process and $Y$ a cadlag process. Let $f\dvtx
{\mathbb R}\times{\mathbb R}\rightarrow{\mathbb R}$ be of
class~$C^{1,3}$. Then
\begin{eqnarray*}
\int_0^t Y\,d^\circ f(V, X) &= &\int_0^t Y\, \partial_v f(V_s,
X_s)\,dV_s + \int_0^t Y\, \partial_x f(V_s, X_s)\,d^\circ X_s
\\
&&{} - \tfrac{1}{12} \int_0^t Y\, \partial^3_{\mathit{xxx}} f(V_s,
X_s)\,d[X,X,X]_s.
\end{eqnarray*}
\end{prop}

We deduce, in particular, that a $C^1$ transformation of a strong zero
cubic variation process is again a strong zero cubic variation process.

We conclude the section by introducing a concept of \textit{definite
integral} via regularization. If processes $X, Y$ are indexed by the
whole real line, a.s. with compact support, we define
\begin{eqnarray}
\label{E1.1def} \int_{\mathbb R}Y\,d^{-} X & := & \lim_{\varepsilon
\to0+} \int_{\mathbb R}Y_s
\frac{X_{s+\varepsilon}-X_s}{\varepsilon} \,  ds,
\\
\label{E1.1Sdef} \int_{\mathbb R}Y_s \, d^{\circ} X_s & := & \lim
_{\varepsilon\to0+} \int_{\mathbb R}Y_s \frac{X_{s+\varepsilon}-X_{s
-\varepsilon} }{2 \varepsilon}\,   ds,
\end{eqnarray}
where the limit is understood in probability. Integration by
parts\break
[Remark~\ref{R1.1}(c)], Proposition~\ref{PEXVQ7} and the chain rule
property (Proposition~\ref{Crule}) can all  be immediately adapted to
these definite integrals.

\section{The PDE operator $L$}\label{s4}

Let $\sigma$, $b\in C^0(\mathbb{R})$ be such that $\sigma>0$. Without
loss of generality, we will suppose that $b(0) = 0$.

We consider a formal PDE operator of the following type:
\begin{equation}\label{E2.1}
Lg = \frac{\sigma^2}{2} g'' + b' g'.
\end{equation}
If $b$ is of class $C^1$, so that $b'$ is continuous, we will say that
$L$ is a \textit{classical} PDE operator.

For a given mollifier $\Phi$, we denote
\[
 \sigma_n^2 := (\sigma^2 \wedge n) \ast\Phi_n,\qquad
  b_n := \bigl(-n \wedge( b \vee n)\bigr) \ast\Phi_n.
\]
We then consider
\begin{eqnarray}\label{E2.2}
L_n g &= & \frac{\sigma_n^2}{2} g'' + b_n' g'\qquad\mbox{for }  g \in
C^2({\mathbb R}), \nonumber
\\[-8pt]
\\[-8pt]
\nonumber {\mathcal  L}_n u &=& \partial_t u + L_n u\qquad\mbox{for } u
\in C^{1,2}([0,T[\times{\mathbb R}),
\end{eqnarray}
where $L_n$ acts on $x$. A priori, $\sigma_n^2$, $b_n$ and the operator
$L_n$ depend on the mollifier~$\Phi$.

Previous definitions are slightly different from those in papers
\cite{frw1,frw2}, but a considerable part of the analysis of $L$ and
the study of the martingale problem can be adapted. In those papers,
there was only regularization but no truncation; here, truncation is
used to study the associated parabolic equations.

\begin{defi} \label{DC1gen}
A function $f\in C^1(\mathbb{R})$ is said to be a
\textit{$C^1$-generalized solution} to
\begin{equation}\label{E2.3}
Lf = \dot{\ell},
\end{equation}
where $\dot{\ell}\in C^0$ if for any mollifier $\Phi$, there are
sequences $(f_n)$ in $C^2$ and  $(\dot{\ell}_n)$ in $C^0$ such that
\begin{equation}\label{E2.4}
L_n f_n = \dot{\ell}_n,\qquad  f_n \to f \mbox{ in } C^1,\qquad
\dot{\ell}_n \to\dot{\ell} \mbox{ in } C^0.
\end{equation}
\end{defi}

\begin{prop}\label{T24}
There is a solution $h\in C^1$ to $Lh=0$ such that $h'(x)\neq0$ for
every $x\in
\mathbb{R}$ if and only if
\[
\Sigma(x) := \lim_{n\to\infty} 2\int_0^x \frac{b_n'}{\sigma
_n^2}(y)\,dy
\]
exists in $C^0$, independently of the mollifier. Moreover, in this
case, any solution $f$ to $Lf=0$ fulfills
\begin{equation}\label{E2.5}
f'(x) = e^{-\Sigma(x)} f'(0).
\end{equation}
\end{prop}

\begin{pf}
This result follows in a very similar way to the proof of
Proposition~2.3 in~\cite{frw1}---first at the level of regularization
and then passing to the limit.
\end{pf}

For the remainder of this paper, we will suppose the existence of this
function~$\Sigma$. We will consider $h \in C^1$ such that
\begin{equation} \label{Efh}
h'(x) := \exp(-\Sigma(x)),\qquad h(0) = 0.
\end{equation}
In particular, $h'(0) = 1 $ holds. Even though  we discuss the general
case with related nonexplosion conditions in \cite{frw1}, here, in
order to ensure conservativeness, we suppose that
\begin{eqnarray} \label{Nonexplo}
\int_{- \infty}^0 e^{- \Sigma(x)}\, dx &=& \int^{ \infty}_0 e^{-
\Sigma(x)}\, dx = +\infty,
\nonumber
\\[-8pt]
\\[-8pt]
\nonumber \int_{- \infty}^0 \frac{e^{ \Sigma(x)}}{\sigma^2}\, dx &=&
\int^{ \infty}_0 \frac{e^{ \Sigma(x)}}{\sigma^2}\, dx = + \infty.
\end{eqnarray}
Previous assumptions are of course satisfied if $\sigma$ is lower
bounded by a positive constant and $b$ is constant outside a compact
interval.

Condition (\ref{Nonexplo}) implies that the image set of $h$ is
${\mathbb R}$.

\begin{rem} \label{R2.5} Proposition~\ref{T24} implies uniqueness of
the problem
\begin{equation}\label{E2.9}
Lf = \dot{\ell}, \qquad  f\in C^1,\qquad   f(0) = x_0,\qquad   f'(0) =
x_1
\end{equation}
for every $\dot{\ell} \in C^0$, $x_0,x_1\in\mathbb{R}$.
\end{rem}

\begin{rem} \label{R2.6}
We present four important examples where $\Sigma$ exists:
\begin{longlist}[(a)]
\item[(a)] If $ b(x) = \alpha ( \frac{\sigma^2 (x) }{2} - \frac
{\sigma^2 (0)}{2}  ) $ for some $\alpha\in\,]0,1]$, then
\[
\Sigma(x) = \alpha\log \biggl(\frac{\sigma^{2}(x)}{\sigma ^{2}(0)}
\biggr)
\]
and
\[
h'(x) = \frac{\sigma^{2\alpha} (0)}{\sigma^{2\alpha}(x)}.
\]
If $\alpha= 1$, the operator $L$ can be formally expressed in divergence
form as $L f = (\frac{\sigma^2}{2} f')'.$

\item[(b)] Suppose that $b$ is locally of bounded variation. We then
get
\[
\int_0^x \frac{b_n'}{\sigma_n^2} (y)\,dy = \int_0^x \frac
{db_n(y)}{\sigma_n^2(y)} \to\int_0^x \frac{db}{\sigma^2}
\]
since $db_n \to db$ in the weak-$\ast$ topology and $\frac{1}{\sigma
^2}$ is continuous.

\item[(c)] If $\sigma$ has bounded variation, then we have
\[
\Sigma(x) = -2\int_0^x b \,  d\biggl(\frac{1}{\sigma^2}\biggr) +
\frac{2b}{\sigma^2}(x) - \frac{2b}{\sigma^2}(0) .
\]
In particular, this example contains the case where $\sigma= 1$ for
any~$b$.

\item[(d)] Suppose that $\sigma$ is locally H\"older continuous with
parameter $\gamma$ and that $b$ is locally H\"older continuous with
parameter $\beta$ such that $\beta+ \gamma> 1$. Since $\sigma$ is
locally bounded, $\sigma^2$ is also locally H\" older continuous with
parameter $\gamma$. Proposition~\ref{PYoungreg} implies that
$\sigma_n^2 \rightarrow\sigma^2$ in $C^{\gamma'}$ and $b_n \rightarrow
b$ in $C^{\beta'}$ for every $\gamma' < \gamma$ and $\beta' < \beta$.
Since $\sigma$ is strictly positive on each compact, $\frac{1}{\sigma
_n^2} \rightarrow\frac{1}{\sigma^2}$ in $C^{\gamma'}$. By
Remark~\ref{REXVQ7}, $\Sigma$~is well defined and locally H\" older
continuous with parameter $\beta'$.
\end{longlist}
\end{rem}

Again, the following lemma can be proven at the level of
regularizations; see also Lemma 2.6 in \cite{frw1}.

\begin{lemma}\label{T27}
The unique solution to problem (\ref{E2.9}) is given by
\begin{eqnarray*}
f(0) & = & x_0,
\\
f'(x) & = & h'(x)  \biggl( 2\int_0^x \frac{\dot{\ell}(y)}{(\sigma
^2h')(y)}\, dy + x_1  \biggr).
\end{eqnarray*}
\end{lemma}

\begin{rem} \label{R2.3}
If $b'\in C^0({\mathbb R})$ and $f\in
C^2({\mathbb R})$ is a classical solution to $Lf = \dot{\ell}$, then
$f$ is clearly also a $C^1$-generalized solution.
\end{rem}

\begin{rem} \label{R2.9bis}
Given $\ell\in C^1$, we denote by $T \ell$ the unique $C^1$-generalized
solution $f$ to problem (\ref{E2.9}) with $\dot\ell= \ell'$, $x_0 = 0,
x_1 =0 $. The unique solution to the general problem (\ref{E2.9}) is
given by
\[
f = x_0 + x_1 h + T \ell.
\]
We write $T^{x_1} \ell= T \ell+ x_1 h$, that is, the solution with $x_0
= 0.$
\end{rem}

\begin{rem} \label{R2.9}
Let $f \in C^1$. There is at most one $\dot{\ell} \in C^0 $ such that
$Lf = \dot{\ell}$. In fact, to see this, it is enough to suppose that
$f = 0$. Lemma \ref{T27} implies that
\[
2 \int_0^x \frac{\dot{\ell}}{\sigma^2 h'}(y)\,dy \equiv0.
\]
Consequently, $\dot{\ell}$ is forced to be zero.
\end{rem}

This consideration allows us to define without ambiguity $L\dvtx
{\mathcal {D}}_L$ $\to C^0$, where ${\mathcal {D}}_L$ is the set of all
$f\in C^1(\mathbb{R})$ which are $C^1$-generalized solution to $Lf =
\dot{\ell}$ for some $\dot{\ell}\in C^0$. In particular, $T
\ell\in{\mathcal {D}}_L$.

A direct consequence of Lemma \ref{T27} is the following
useful result.

\begin{lemma}\label{L28bis}
${\mathcal {D}}_L $ is the set of $f \in C^1$ such that there exists
$\psi\in C^1$ with $f' = e^{- \Sigma} \psi.$
\end{lemma}

In particular, it gives us the following density proposition.

\begin{prop}\label{T29}
${\mathcal {D}}_L$ is dense in $C^1$.
\end{prop}

\begin{pf}
It is enough to show that every $C^2$-function is the $C^1$-limit of a
sequence of functions in ${\mathcal {D}}_L$. Let $(\psi_n)$ be a
sequence in $C^1$ converging to $f' e^ \Sigma$ in~$C^0$. It follows
that
\[
f_n (x) := f (0) + \int_0^x e^ {- \Sigma}(y) \psi_n (y)\,dy,\qquad   x
\in{\mathbb R},
\]
converges to $f \in C^1$ and $f_n \in{\mathcal {D}}_L$.
\end{pf}

We must now discuss technical aspects of the way $L$ and its domain
${\mathcal {D}}_L$ are transformed by $h$. We recall that $Lh=0$ and
that $h'$ is strictly positive. Condition~(\ref{Nonexplo}) implies that
the image set of $h$ is ${\mathbb R}$.

Let $L^0$ be the classical PDE operator
\begin{equation}\label{E2.12}
L^0 \phi= \frac{\tilde\sigma_h^2}{2} \phi'',\qquad \phi\in C^2,
\end{equation}
where
\[
\tilde\sigma_h (y) = (\tilde\sigma h')(h^{-1}(y)),\qquad y\in{\mathbb
R}.
\]
$L^0$ is a classical PDE map; however, we can also consider it at the
formal level and introduce ${\mathcal  D}_{L^0}$.

\begin{prop}\label{T212}
\textup{(a)} $h^2\in{\mathcal {D}}_L$, $Lh^2 = h'^2\sigma^2$.
{\smallskipamount=0pt
\begin{longlist}[(b)]
\item[(b)] ${\mathcal {D}}_{L^0} = C^2$.

\item[(c)] $\phi\in{\mathcal {D}}_{L^0}$ holds if and only if
$\phi\circ h \in{\mathcal {D}}_L$. Moreover, we have
\begin{equation}\label{E2.13}
L(\phi\circ h) = (L^0\phi)\circ h
\end{equation}
for every $\phi\in C^2$.
\end{longlist}}
\end{prop}

\begin{pf}
This follows similarly as for Proposition~2.13 of \cite{frw1}.
\end{pf}

We will now discuss another operator related to $L$. Given a function
$f$, we need to provide a suitable definition of $f \mapsto\int_0^x
Lf(y)\,dy$, that is, some primitive of $Lf$.
\begin{itemize}
\item One possibility is to define that map, through previous
expression,
for $f\in{\mathcal  D}_L$.

\item Otherwise, we try to define it as linear map on $C^2$. For this,
first suppose that $b'$ is continuous. Then integrating by parts, we
obtain
\begin{equation} \label{E213}
\int_0^x Lf(y)\,dy = \int_0^x  \biggl(\frac{\sigma^2}{2} - b
\biggr)f''(y)\,dy + (bf')(x) - (bf')(0).
\end{equation}
We remark that the right-hand side of this expression makes sense for any
$f \in C^2$ and continuous $b$. We will thus define $\hat L \dvtx C^2
\rightarrow C_0^0$ as follows:
\begin{equation} \label{E213a}
\hat L f : = \int_0^x  \biggl(\frac{\sigma^2}{2} - b  \biggr)
f''(y)\,dy + (bf')(x) - (bf')(0).
\end{equation}
\end{itemize}
One may ask if, in the general case, the two definitions $f
\rightarrow\int_0^x Lf(y)\,dy$ on ${\mathcal  D}_L$ and $\hat L$ on
$C^2$ are compatible. We will later see that under Assumption~\ref{A0}, 
this will be the case. However, in general,
${\mathcal D}_L \cap C^2$ may be empty.

Thus far, we have learned how to eliminate the first-order term in a
formal PDE operator through the transformation $h$ introduced at
(\ref{Efh}); when $L$ is classical, this was performed by Zvonkin (see
\cite{z}). We would now like to introduce a transformation which puts
the PDE operator in a divergence form.

Let $L$ be a PDE operator which is formally of type (\ref{E2.1}):
\[
L g = \frac{\sigma^2}{2} g'' + b' g'.
\]

We consider a function of class $C^1$, namely $k\dvtx
\mathbb{R}\to\mathbb {R}$ such that
\begin{equation} \label{Efk}
k(0) = 0\quad\mbox{and}\quad  k'(x) = \sigma^{-2}(x)\exp(\Sigma(x)).
\end{equation}
According to assumptions (\ref{Nonexplo}), $k$ is bijective on
${\mathbb R}$.

\begin{rem} \label{R2.19}
If there is no drift term, that is, $b = 0$, then we have
$k'(x)=\sigma^{-2}(x)$.
\end{rem}

\begin{lemma}\label{T216}
We consider the formal PDE operator
given by
\begin{equation}\label{E2.21}
L^1 g = \frac{{\bar\sigma}_k^2}{2} g'' +
\biggl(\frac{{\bar\sigma}^2_k}{2} \biggr)' g' =
\biggl(\frac{{\bar\sigma}^{2}_k}{2} g'\biggr)',
\end{equation}
where
\[
\bar\sigma_k (z) = (\sigma k') \circ k^{-1}(z), \qquad  z \in{\mathbb
R}.
\]
Then:
\begin{longlist}[(ii)]
\item[(i)] $g\in{\mathcal {D}}_{L^1}$ if and only if $ g\circ
k\in{\mathcal {D}}_L$;

\item[(ii)] for every $g\in{\mathcal {D}}_{L^1}$, we have $L^1 g =
L(g\circ k)\circ k^{-1}$.
\end{longlist}
\end{lemma}

\begin{pf}
It is practically the same as in Lemma 2.16 of \cite{frw1}.
\end{pf}

  We now give a lemma whose proof can be easily established
by investigation. Suppose that $L$ is a classical PDE operator. Then
${\mathcal  L}= \partial_t + L$ is well defined for $C^{1,2}([0,T[
\times{\mathbb R})$ functions where $L$ acts on the second variable.
Given a function $\varphi\in C([0,T] \times{\mathbb R})$, we will
hereafter set $\tilde{\varphi}\dvtx [0,T] \times{\mathbb
R}\longrightarrow{\mathbb R}$ by $\tilde{\varphi}(t,y)
=\varphi(t,h^{-1}(y))$.

\begin{lemma} \label{L28}
Let us suppose that $h\in C^{2}({\mathbb R})$. We set $\sigma_h =\sigma
h'$.

We define the PDE operator ${{\mathcal  L}}^0$ by ${{\mathcal
L}}^0\varphi=
\partial_t \varphi+ L^0 \varphi $, where $L^0$ is a classical operator
acting on the space variable $x$ and
\[
L^0 f =
\frac{\tilde{\sigma}_h^2}{2} f''.
\]
If $f\in C^{1,2}([0,T[ \times{\mathbb R})$ and ${\mathcal  L}f= \gamma$
in the classical sense, then ${{\mathcal  L}}^0 \tilde{f} =
\tilde{\gamma}$.
\end{lemma}

We will now formulate a supplementary assumption which will be useful
when we study singular stochastic differential equations in the proper
sense and not only in the form of a martingale problem.

\renewcommand{\theta}{${\mathcal  A}(\nu_0$)}

\begin{ta}\label{A0}
Let $\nu_0$ be a topological F-space which is a linear topological
subspace of $C^0 ({\mathbb R})$ (or, eventually, an inductive limit of
sub-F-spaces). The $\nu_0$-convergence implies convergence in $C^0$
and, therefore, pointwise convergence.

We say that $L$ fulfills Assumption \textup{\ref{A0}} if the following
conditions hold:

\begin{longlist}[(iii)]
\item[(i)] $C^1 \subset\nu_0$, which is dense. 

\item[(ii)] For every $g \in C^1({\mathbb R})$, the multiplicative
operator $\phi\rightarrow g \phi$ maps $\nu_0$ into itself.

\item[(iii)] Let $T\dvtx  C^1({\mathbb R}) \subset\nu_0 \rightarrow
C^1({\mathbb R})$ as defined in Lemma \ref{T27}, that is, $f = T \ell$
is such that
\begin{eqnarray*}
f(0) & = & 0,
\\
f'(x) & = & e^{-\Sigma}(x)  \biggl( 2\int_0^x \frac{e^{\Sigma}(y)
{\ell}'(y)}{\sigma^2(y)}  \,dy  \biggr).
\end{eqnarray*}
We recall that $f = T \ell$ solves problem $Lf = \ell'$ with $f(0) =
f'(0) = 0$. We suppose that $T$ admits a continuous extension to
$\nu_0$.

\item[(iv)] Let $x_1 \in{\mathbb R}$. For every $f \in C^2$ with $f(0)
= 0$ and $ f'(0) = x_1$ so that ]\mbox{$\hat Lf = \ell$}, we have
$\ell\in\nu_0$ and $T^{x_1} \ell=f$, where $T^{x_1} $ denotes the
continuous extension of $T^{x_1}$ (see Remark~\ref{R2.9bis}) to
$\nu_0$, which exists by~\textup{(iii)}.

\item[(v)]  The set $\hat L C^2$ is dense in $\{\ell\in\nu_0 |\ell (0)
= 0 \}$.
\end{longlist}
\end{ta}

\begin{rem} \label{R2.21a} Let $x_1 \in{\mathbb R}$.
\begin{longlist}[(iii)]
\item[(i)] Remark~\ref{R2.9bis} and point (iii) above together imply
that $T^{x_1} \dvtx C^1({\mathbb R}) \subset\nu_0 \rightarrow
C^1({\mathbb R}) $ extends continuously to $\nu_0$. Moreover,
\[
\{f \in C^2 | f(0) = 0, f'(0) = x_1 \} \subset Im T^{x_1}.
\]

\item[(ii)] Point (iv) above shows that $b \in \nu_0$ and $ T^1 b =
id$, where $id(x) = x$; in fact, $id(0) = 0$, $id'(1) = 1$ and
(\ref{E213a}) implies that $\hat L id = b$.

\item[(iii)] Point (i) above is satisfied if, for instance, the map $T$
is closable as a map from $C^0$ to $C^1$. In that case, $\nu_0$ may be
defined as the domain of the closure of $C^1$, equipped with the graph
topology related to $C^0 \times C^1$.
\end{longlist}
\end{rem}

Below, we give some sufficient conditions for points (iv) and (v) of
the Technical Assumption to be satisfied.

We define by $ C^1_{\nu_0} $ the vector space of functions $f \in C^1$
such that $f' \in\nu_0$. This will be an F-space if equipped with the
following topology. A sequence $(f_n)$ will be said to converge to $f$
in $ C^1_{\nu_0} $ if $f_n (0) \rightarrow f(0)$ and $(f'_n)$ converges
to $f'$ in $\nu_0$. In particular, a sequence converging according to $
C^1_{\nu_0} $ also converges with respect to~$C^1$. On the other hand,
$C^2 \subset C^1_{\nu_0} $ and a sequence converging in $ C^2 $ also
converges with respect to $C^1_{\nu_0}$. Moreover, $C^2$ is dense in
$C^1_{\nu_0}$ because $C^1$ is dense in~$\nu_0$.

\begin{lemma}\label{L2.21}
Suppose that points \textup{(i)} to \textup{(iii)} of the Technical
Assumption are fulfilled. We suppose, moreover, that:
\begin{longlist}[(a)]
\item[(a)] $ h \in C^1_{\nu_0}$.

\item[(b)] For every $f \in C^2$, $f(0) = 0$, $f'(0) = 0$, $\hat Lf =
\ell$, we have $\ell\in\nu_0$ and \mbox{$T \ell=f$}.

\item[(c)] $\hat L\dvtx  C^2 \rightarrow\nu_0$ is well defined and has
a continuous extension to $ C^1_{\nu_0} $, still denoted by $\hat L$,
such that $\hat L h = 0$.

\item[(d)] $Im T \subset C^1_{\nu_0}$.

\item[(e)] $\hat L T $ is the identity map on $ \{ \ell\in\nu_0
|\ell(0) = 0 \} $.
\end{longlist}

Then $T, T^{x_1} $ for every $x_1 \in{\mathbb R}$ are injective and
points \textup{(iv)} and \textup{(v)} of the Technical Assumption are
satisfied.
\end{lemma}

\begin{pf}
The injectivity of $T$ follows from point (e). The injectivity of
$T^{x_1}$ is a consequence of Remark~\ref{R2.9bis}.

We prove point (iv). Point (c) says that $\hat L h = 0$. We set $\hat f
= f - x_1 h$, $f\in C^2$, where $f(0)=0$, $f'(0)=x_1$. Clearly, $ \hat
L \hat f = \hat L f = \ell$ and $\hat f (0) = 0, \hat f' (0) = 0$.
Point (b) implies that $T \ell= \hat f$. Hence, $T^{x_1} \ell=
T\ell+x_1 h=f $ and (iv) is satisfied.

Concerning point (v), let $\ell\in\nu_0$ with $\ell(0) = 0$ and set $f
= T \ell$. Since $f$ belongs to $ C^1_{\nu_0} $ by (c), $f'$ belongs to
$\nu_0$. Point (i) of the technical assumption implies that there
exists a sequence $(f'_n)$ of $C^1$ functions converging to $f'$ in the
$\nu_0$ sense and thus also in $C^0$. Let $(f_n)$ be the sequence of
primitives of $(f'_n)$ (which are of class~$C^2$) such that $f_n(0) =
0$. In particular, we have that $(f_n)$ converges to $f$ in the
$C^1_{\nu_0}$-sense. By (c), there exists $\lambda$ in $\nu_0$ which is
the limit of $\hat L f_n$ in the $\nu_0$-sense. Observe that because
of~(b), $T ( \hat L f_n)=f_n$. On the other hand, $\lim_{n \rightarrow+
\infty} f_n = f$ in $C^1$. Applying $T$ and using (iii) of the
Technical Assumption, we obtain
\[
T \lambda= \lim_{n \rightarrow+ \infty} T ( \hat L f_n) = \lim_{n
\rightarrow+ \infty} f_n = f = T \ell.
\]
The injectivity of $T$ allows us to conclude that $\ell= \lambda$.
\end{pf}

\begin{rem} \label{R2.21aa}
Under the assumptions of Lemma \ref{L2.21}, we have:
\begin{itemize}
\item ${\mathcal  D}_L \subset C^1_{\nu_0}$;

\item $\hat L f = \int_0^x Lf(y)\,dy$, $ f \in{\mathcal  D}_L.$
\end{itemize}

In fact, let $f \in{\mathcal  D}_L$. Without loss of generality, we can
suppose that $f(0) = 0$. Let $x_1 = f'(0)$ and set $\hat f = f + x_1 h$
so that $ \hat f(0) = \hat f'(0) = 0$. Setting $\dot\ell= L \hat f,$
Lemma \ref{T27} implies that $\hat f = T \ell$, where $\ell= \int _0^x
\dot\ell(y)\,dy$. So $\hat f \in Im T \subset C^1_{\nu_0}$. Since $h
\in C^1_{\nu_0}$, it follows that $f \in C^1_{\nu_0}$, by additivity.

On the other hand,
\begin{eqnarray*}
L f &=& L \hat f + x_1 Lh = \hat L f = \dot\ell,
\\
\hat L f &=& \hat L \hat f + x_1 \hat Lh = \hat L T \ell= \ell,
\end{eqnarray*}
by point (e) of Lemma \ref{L2.21}.
\end{rem}

\begin{exemple} \label{4ex}
We provide here a series of four significant examples when Technical
Assumption \textup{\ref{A0}} is verified. We only comment on the points
which are not easy to verify.
\begin{longlist}[(iii)]
\item[(i)] The first example is simple. It concerns the case when the
drift $b'$ is continuous. This problem, to be studied later,
corresponds to an ordinary SDE where
\[
\nu_0 = C^1,\qquad C^1_{\nu_0} = C^2,\qquad \hat L f =
\int_0^{\bolds{\cdot}} Lf(y)\,dy.
\]

\item[(ii)] $L$ is \textit{close to divergence type}, that is, $b =
\frac{\sigma^2 - \sigma^2(0)}{2} + \beta$ and where $\beta$ is a locally bounded
variation function vanishing at zero. The operator is of divergence type with an
additional Radon measure term, that is, we have $\Sigma= \ln\sigma^2 +
2 \int_0^x \frac{d\beta}{\sigma ^2} $. In this case, we have $\nu_0 =
C^0.$ Points (i) and (ii) of the Technical Assumption are trivial.

We have, in fact,
\[
h'(x) = e^{-\Sigma} = \frac{1}{\sigma^2(x)} \exp \biggl(-2 \int_0^x
\frac{d\beta}{\sigma^2} \biggr).
\]
$T$ defined at point (iii) of the Technical Assumption is such that $T
\ell= f$, where $f(0) = 0$ and
\begin{equation}\label{E2.26}
f'(x) = \frac{2 \sigma^2 (0)}{\sigma^2(x)} \exp \biggl(-2 \int_0^x
\frac{d\beta}{\sigma^2} \biggr) \int_0^x \ell'(y)\exp \biggl(2\int_0^y
\frac{d\beta}{\sigma ^2} \biggr)\,dy.
\end{equation}
Consequently, the extension of $T$ to $\nu_0 = C^0$, always still denoted by
the same letter $T$, is given by $f = T \ell$ with $f(0) = 0$ and
\begin{eqnarray}\label{E2.27}
\qquad f'(x) & = & \frac{2}{\sigma^2(x)}  \biggl\{ \ell(x) - 2 \exp \biggl(-2
\int_0^x \frac{d\beta}{\sigma^2}
 \biggr)
 \nonumber
 \\[-8pt]
 \\[-8pt]
 \nonumber
\qquad &&\hspace*{11.2mm} {}\times \biggl( \ell(0) +  \int_0^x \ell(y) \exp
\biggl(2\int_0^y\frac{d\beta}{\sigma ^2} \biggr) \frac{1}{\sigma^2(y)}
\,d\beta(y) \biggr) \biggr\}.
\end{eqnarray}
Points (iv) and (v) are seen to be satisfied via Lemma \ref{L2.21}. We
have $C^1_{\nu_0} = C^1$. Point (a) is obvious since $h' \in C^0$ and
so $h \in C^1_{\nu_0}$. Let $f \in C^2$. Using Lebesgue--Stieltjes
calculus, we can easily show that
\begin{equation} \label{E2.27a}
\ell(x) = \hat L f (x) = \frac{\sigma^2 (x)}{2} f'(x) - \frac {\sigma^2
(0)}{2} f'(0) + \int_0^x f' \,d\beta.
\end{equation}
This shows that $\ell\in C^0 = \nu_0$ and therefore the first part of
(b). We remark that we can, in fact, consider $\hat L\dvtx C^2
\rightarrow\nu_0$ because
\[
\hat L f = \hat L (f-x_1 h) + x_1 \hat L h = \hat L (f-x_1 h) \in\nu
_0.
\]
The expression of $\hat L f$ extends continuously to $f \in C^1$, which
yields the first part of point~(c). Moreover, inserting the expression
for $h'$ into $f'$ in (\ref{E2.27a}), one shows that $\hat L h = 0$.

Suppose, now, that in expression (\ref{E2.27a}), $f \in C^2$, $f(0) =
0$, $f'(0) = 0$. A simple investigation shows that $T \ell= f$, so the
second part of point (b) is fulfilled; point~(d) is also clear because
of (\ref{E2.27}). Finally point (d) holds because one can prove by
inspection that $\hat L T $ is the identity on $C^0_0$.

\item[(iii)] We recall the notation $D^\gamma({\mathbb R})$ which
indicates the topological vector space of locally H\"older continuous
functions defined on ${\mathbb R}$ with parameter $\alpha> \gamma$. We
recall that $D^\gamma({\mathbb R})$ is a vector algebra.

Suppose that $\sigma\in D^{{1}/{2}}$ and $b \in C^{{1}/{2}}$ (or
$\sigma\in C^{{1}/{2}}$ and $b \in D^{{1}/{2}}$). Remark~\ref{R2.6}(d)
implies that $\Sigma$ also belongs to $D^{{1}/{2}}$. We set $\nu_0 =
D^{{1}/{2}}$.

Technical Assumption \ref{A0} is verified for the following reasons.

Since $\Sigma\in D^{{1}/{2}}$, $ h' = e^{-\Sigma}$ belongs to the same
space.

Point (i) follows because of Proposition~\ref{PYoungreg} and point (ii)
follows because $D^{{1}/{2}}$ is an algebra. Corollary \ref{CYoung}
yields that for every $\ell\in D^{{1}/{2}}$, the function
\begin{equation} \label{EImT}
f'(x) = e^{-\Sigma(x)} \int_0^x 2 \frac{e^\Sigma}{\sigma^2}(y)\,d^{(y)}
\ell(y)
\end{equation}
is well defined and belongs to $D^{{1}/{2}}$. This shows that $T$ can
be continuously extended to $\nu_0$ and point (iii) is established.

Concerning points (iv) and (v), we again use Lemma \ref{L2.21}. We
observe that
\[
C^1_{\nu_0} = \{ f \in C^1 | f' \in D^{{1}/{2}} \}.
\]
Point (a) is obvious since $h' = e^{- \Sigma} \in D^{{1}/{2}}$. Let $f
\in C^2$. Considering $b$ as a deterministic process and
recalling the definition of $\hat L$ as in (\ref{E213a}), integration
by parts in Remark~\ref{R1.1}(c) and Proposition~\ref{PEXVQ7} together
imply that
\begin{eqnarray} \label{E1.1b}
\ell(x) & = & \int_0^x \frac{\sigma^2}{2} \,d^{0} f' + \int_0^x f'
\,d^\circ b,
\\
\ell(x) & = & \int_0^x \frac{\sigma^2}{2} \,d^{(y)} f' + \int_0^x f'
\,d^{(y)} b.
\end{eqnarray}
The first part of point (b) follows because of
Proposition~\ref{PCRuleY}. Of course, the previous expression can be
extended to $f \in C^1_{\nu_0}$ and this shows the first part of
point~(c).

Showing that the second part of point (c) of Lemma \ref{L2.21} holds
consists of  verifying that $\hat L h = 0$. Substituting $h' =
e^{-\Sigma}$ into the previous expression, through
Proposition~\ref{PCRuleY}, we obtain
\[
\ell(x) = - \int_0^x \frac{\sigma^2}{2} e^{-\Sigma} \,d^{(y)} \Sigma+
\int_0^x e^{-\Sigma}\, d^{(y)} b = 0.
\]
Concerning the second part of point (b), let $f \in C^2$ so that $f(0)
= f'(0) = 0$. We want to show that $\varphi= T \ell$ coincides with
$f$.

Since $\varphi(0) = 0$, it remains to check that $ \varphi' = f'$. We
recall that
\[
\varphi'(x) = e^{-\Sigma}(x)  \biggl( 2\int_0^x \frac{e^\Sigma}
{\sigma^2}(y)\,d^{(y)} \ell(y)  \biggr).
\]

Twice applying the chain rule of Proposition~\ref{PCRuleY} and using
(\ref{E1.1b}), the fact that
\[
e^\Sigma(x) = \int_0^x e^\Sigma\frac{2 \,d^{(y)} b}{\sigma^2} + 1
\]
and integration by parts, we obtain
\begin{eqnarray*}
\varphi'(x) &=& e^{-\Sigma}(x)  \biggl\{ \int_0^x e^\Sigma \,d^{0} f' +
\int_0^x 2 \frac{e^\Sigma}{\sigma^2} f' \,d^{(y)} b  \biggr\}
\\
&=& e^{-\Sigma}(x)  \biggl\{ \int_0^x e^\Sigma \,d^{0} f' +
\int_0^x f' \,d^{(y)} e^\Sigma \biggr\}
\\
&=& e^{-\Sigma}(x)  \biggl\{ \int_0^x e^\Sigma \,d^{0} f' +
\int_0^x f' \,d^{0} e^\Sigma \biggr\}
\\
&=& e^{-\Sigma}(x)  \{ (f' e^\Sigma)(x) - (f' e^\Sigma)(0)
 \} \\
&=& f'(x) .
\end{eqnarray*}
Point (b) is therefore completely established.

Point (d) follows because in (\ref{EImT}), when $\ell\in\nu_0$, it
follows that  $f' \in\nu_0$.

Clearly, as for the previous example, $Im T \subset C^1_{\nu_0}$. It
remains to show that $\hat L T$ is the identity map $\{f \in D^{1/2} |
f(0) = 0 \} $.

For this, we first remark that
\begin{equation} \label{EYVerif1}
\hat L f (x) = \int_0^x \frac{\sigma^2}{2} e^{- \Sigma}\, d^{(y)}
 ( f' e^{ \Sigma}  ).
\end{equation}
In fact, by Proposition~\ref{PEXVQ7} and integration by parts contained
in Remark~\ref{R1.1}(c), we obtain
\[
f' (x) e^{\Sigma(x)} = f'(0) + \int_0^x e^{\Sigma}\, d^{(y)} f' +
\int_0^x f'\, d^{(y)} e^{ \Sigma}.
\]
By the chain rule of Proposition~\ref{PCRuleY}, we obtain the
right-hand side of (\ref{EYVerif1}).

At this point, by definition, if $f = T\ell$, we have
\[
f' (x) e^{\Sigma(x)} = \int_0^x 2 \frac{e^{ \Sigma}}{\sigma^2}
\,d^{(y)} \ell.
\]
Therefore, (\ref{EYVerif1}) and Proposition~\ref{PCRuleY} allow us to
conclude that
\[
\hat L f (x) = \int_0^x \frac{\sigma^2}{2} e^{- \Sigma} 2 \frac
{e^{\Sigma}}{\sigma^2} \,d^{(y)} \ell= \ell(x)-\ell(0) .
\]

\item[(iv)] Suppose $b$ is locally with bounded variation. Then the
Technical Assumption is satisfied for $\nu_0 = \mathit{BV}$, where
$\mathit{BV}$ is the space of continuous real functions, locally with
bounded variation $v$, equipped with the following topology. A sequence
$(v_n)$ in $\mathit{BV}$ converges to $v$ if
\begin{eqnarray*}
v_n(0) &\rightarrow & v(0),
\\
dv_n &\rightarrow &dv \qquad \mbox{in  the  weak- $*$ topology}.
\end{eqnarray*}
The arguments for proving that the Technical Assumption is satisfied
are similar, but easier, than those for the previous point. Young-type
calculus is replaced by classical Lebesgue--Stieltjes calculus.
\end{longlist}
\end{exemple}

\section{Martingale problem}\label{s5}

In this section, we consider a PDE operator satisfying the same
properties as in previous section, that is,
\begin{equation} \label{E3.1}
Lg = \frac{\sigma^2}{2} g'' + b'g',
\end{equation}
where $\sigma>0$ and $b$ are continuous. In particular, we assume that
\begin{equation} \label{E3.2}
\Sigma(x) = \lim_{n\to\infty} 2\int_0^x \frac{b_n'}{\sigma _n^2}(y)\,dy
\end{equation}
exists in $C^0$, independently of  the chosen mollifier. Then $h$
defined by $h'(x):=\exp(-\Sigma(x))$ and $h(0)=0$ is a solution to
$Lh=0$ with $h'\neq0$.

Here, we aim to introduce different notions of martingale problem,
trying, when possible, to also clarify the classical notion. For the
next two definitions, we consider the following convention. Let
$(\Omega, {\mathcal  F}, P)$ equipped with a filtration $({\mathcal
F}_t)_{t\geq0}$ fulfill the \textit{usual conditions}; see, for instance, \cite{ks}, Definition
2.25, Chapter 1.

\begin{defi}\label{Dmp}
A process $X$ is said to solve \textit{the martingale problem} related
to $L$ (with respect to the aforementioned filtered probability space)
with initial condition $X_0=x_0$, $x_0\in\mathbb{R}$, if
\[
f(X_t) - f(x_0) - \int_0^t Lf(X_s)\,ds
\]
is an $({\mathcal  F}_t)_{t\geq s}$-local martingale for $f\in{\mathcal
{D}}_L$ and $X_0=x_0$.

More generally, for $s\ge0$, $x\in\mathbb{R}$, we say that $(X_t^{s,x},
t\ge0)$ solves the martingale problem related to $L$ with initial value
$x$ at time $s$ if for every $f\in{\mathcal {D}}_L$,
\[
f(X_t^{s,x}) - f(x) - \int_s^t Lf(X_r^{s,x})\,dr,\qquad   t\ge s,
\]
is an $({\mathcal  F}_t)_{t \geq s}$-local martingale.
\end{defi}

We remark that $X^{s,x}$ solves the martingale problem at time $s$ if and
only if $X_t:= X_{t+s}^{s,x}$ solves the martingale problem
at time $0$.

\begin{defi}\label{DSmp} Let $(W_t)$ be an $({\mathcal  F}_t)$-classical
Wiener process. An\break \mbox{$({\mathcal  F}_t)$-}progressively
measurable process $X = (X_t)$ is said to solve \textit{the sharp
martingale problem} related to $L$ (on the given filtered probability
space) with initial condition $X_0=x_0$, $x_0\in\mathbb{R}$, if
\[
f(X_t) - f(x_0) - \int_0^t Lf(X_r)\,dr = \int_0^t f'(X_r) \sigma
(X_r)\,dW_r
\]
for every $f\in{\mathcal {D}}_L$.

More generally, for $s\ge0$, $x\in\mathbb{R}$, we say that $(X_t^{s,x},
t\ge s)$ solves the sharp martingale problem related to $L$ with
initial value $x$ at time $s$ if for every $f\in{\mathcal {D}}_L$,
\[
f(X_t^{s,x}) - f(x) - \int_s^t Lf(X_r^{s,x})\,dr = \int_s^t f'(X^{s,x}
_r) \sigma(X^{s,x} _r)\,dW_r, \qquad  t\ge s.
\]
\end{defi}

\begin{rem} \label{R3.1a}
Let $(W_t)$ be an $({\mathcal  F}_t)$-Wiener process. If $b'$ is
continuous, then a process $X$ solves the (corresponding) sharp
martingale problem with respect to $L$ if and only if it is a classical
solution of the SDE
\[
X_t = x_0 + \int_0^t b'(X_r)\,dr + \int_0^t \sigma(X_r)\,dW_r.
\]
For this, a simple application of the classical It\^o formula gives the result.
\end{rem}

\begin{rem} \label{R3.2}
(i) In general, $f(x)=x$ does not belong to ${\mathcal {D}}_L$,
otherwise a solution to the martingale problem with respect to $L$
would be a semimartingale. According to Remark~\ref{RSem}, this is
generally not the case. In \cite{frw2}, we gave necessary and
sufficient conditions on $b$ so that $X$ is a semimartingale.

(ii) Given a solution $X$ to the martingale problem related to $L$, we
are interested in the operators
\[
{\mathcal {A}} \dvtx  {\mathcal {D}}_L \to{\mathcal {C}},\qquad
\mbox{given by } {\mathcal {A}}(f) = \int_0^{\bolds{\cdot}}
Lf(X_s)\,ds,
\]
and
\[
A \dvtx  C^1 \to{\mathcal {C}},\qquad\mbox{given by }  A(\ell) =
\int_0^{\bolds{\cdot}} \ell'(X_s)\,ds,
\]
where ${\mathcal  C}$ is the vector algebra of continuous processes.

We may ask whether $\mathcal {A}$ and $A$ are closable in $C^1$ and
$C^0$, respectively. We will see that $\mathcal {A}$ admits a
continuous extension to $C^1$. However, $A$ can be extended
continuously to some topological vector subspace $\nu_0$ of $C^0$,
where $\nu_0$ includes the drift, only when Assumption~\ref{A0} is
satisfied.
\end{rem}

Similarly, as in the case of classical stochastic differential
equations, it is possible to distinguish two types of existence and
uniqueness for the martingale problem. Even if we could treat initial
conditions which are random ${\mathcal F}_0$-measurable solutions,
here we will only discuss deterministic ones. We will denote by $\mathit{MP}(L,
x_0)$ [resp. $\mathit{MP}(L,x_0)$] the martingale problem (resp. sharp martingale problem) related to $L$ with initial condition
$x_0.$ The notions will only be formulated with respect to the initial
condition at time~0.

\begin{defi}[(\textit{Strong existence})]\label{D112}
We will say that $\mathit{SMP}(L,x_0)$ admits \textit{strong existence} if the
following holds. Given any probability space $(\Omega, {\mathcal  F},
P)$, a filtration $({\mathcal  F}_t)_{t \ge0}$ and an $({\mathcal
F}_t)_{t \ge0}$-Brownian motion $(W_t)_{t \ge0}, x_0 \in{\mathbb R}$,
there is a process $(X_t)_{t \ge0}$ which solves the sharp martingale  problem
with respect to $L$ and initial condition $x_0$.
\end{defi}

\begin{defi}[(\textit{Pathwise uniqueness})]\label{D113}
We will say that $\mathit{SMP}(L,x_0)$ admits \textit{pathwise uniqueness} if the
following property is fulfilled.

Let $(\Omega, {\mathcal  F}, P)$ be a probability space with filtration
$({\mathcal  F}_t)_{t \ge0}$ and\break \mbox{$({\mathcal  F}_t)_{t \ge0}$-}Brownian
motion $(W_t)_{t \ge0}$. If two processes $X, \tilde X$ are two
solutions of the sharp martingale  problem with respect to $L$ and
$x_0$,
such that $X_0 = \tilde X_0$ a.s., then $X$ and $\tilde X$
coincide.
\end{defi}

\begin{defi}[(\textit{Existence in law} or \textit{weak existence})] \label{D114}
We will say that $\mathit{MP}(L;x_0)$ \textit{admits weak existence} if there is a
probability space $(\Omega, {\mathcal  F}, P)$, a~filtration
$({\mathcal F}_t)_{t \ge0}$ and a process $(X_t)_{t \ge0}$ which is a
solution of the corresponding martingale problem.

We say that $ \mathit{MP}(L)$ admits weak existence if $\mathit{MP}(L;x_0)$ admits weak existence
for every $x_0$.
\end{defi}

\begin{defi} [(\textit{Uniqueness in law})]\label{D115}
We say that $\mathit{MP}(L;x_0)$ has a \textit{unique solution in law} if the
following holds. We consider an arbitrary probability space $(\Omega,
{\mathcal  F}, P)$ with a filtration $({\mathcal  F}_t)_{t \ge0}$ and a
solution $X$ of the corresponding martingale problem. We
also consider another probability space $(\tilde\Omega, \tilde{\mathcal
F}, \tilde P)$ equipped with another filtration $(\tilde{\mathcal
F}_t)_{t \ge0}$ and a solution $\tilde X$. We suppose that $X_0 = x_0$, \mbox{$P$-a.s.} and
$\tilde X_0 = x_0$, $\tilde P$-a.s. Then $X$ and $\tilde X$ must have
the same law as a r.v.'s with values in $E = C({\mathbb R}_+)$ (or
$C[0,T]$).

\end{defi}


\begin{rem} \label{RSclass}
Let us suppose $b'$ to be a continuous function. We do not suppose
$\sigma$ to be strictly positive (only continuous).
\begin{longlist}[(ii)]
\item[(i)] The $\mathit{SMP}(L,x_0)$ then admits strong existence and pathwise
uniqueness if the corresponding classical SDE
\[
X_t = x_0 + \int_0^t \sigma(X_s)\,dW_s + \int_0^t b'(X_s)\,ds
\]
admits strong existence and pathwise uniqueness. In this case,
${\mathcal  D}_L = C^2$ and to establish this, it is enough to use the
classical It\^o formula.

\item[(ii)] It is well known (see \cite{ks,SV}) that weak existence
(resp., uniqueness in law) of the martingale problem is equivalent to
weak existence (resp., uniqueness in law) of the corresponding SDE.
\end{longlist}
\end{rem}

For the rest of the section let $s \in[0,T]$, $x_0 \in{\mathbb R}. $
Moreover, let $(\Omega, ({\mathcal  F}_t), P) $ be a fixed filtered
probability space fulfilling the usual conditions.

The first result concerning solutions to the martingale problem related
to $L$ is the following.

\begin{prop} \label{P3.3}
Let $ y_0 = h(x_0).$
\begin{longlist}[(ii)]
\item[(i)] A process $X$ solves the martingale problem related to $L$
with initial condition $x$ at time $s$ if and only if $Y=h(X)$ is a
local martingale which solves, on the same probability space,
\begin{equation} \label{E3.3}
Y_t = y_0 + \int_s^t \tilde\sigma_h (Y_s)\,dW_s,
\end{equation}
where $\tilde\sigma_h(y) = (\sigma h')(h^{-1}(y))$ and where $(W_t)$ is
an $({\mathcal  F}_t)$-classical Brownian motion.

\item[(ii)] Let $(W_t)$ be an $({\mathcal  F}_t)$-classical Brownian
motion. If $Y$ is a solution to equation~(\ref{E3.3}), then $X =
h^{-1}(Y)$ is a solution to the sharp martingale problem with respect
to $L$ with initial condition $x$ at time $s$.
\end{longlist}
\end{prop}

\begin{rem} \label{R3.4}
Let $X$ be a solution to the martingale problem with respect to $L$ and
set $Y = h(X)$ as in point~(i) above. Since $Y$ is a local martingale,
we know from Remark~\ref{R81}(iv) that $X=h^{-1}(Y)$ is an $({\mathcal
F}_t)$-Dirichlet process with martingale part
\[
M_t^X = \int_0^t (h^{-1})'(Y_s)\,dY_s.
\]
In particular, $X$ is a finite quadratic variation process with
\[
[X,X] =[M^X,M^X]_t = \int_0^t \sigma^2(X_s)\,ds.
\]
\end{rem}

\begin{pf*}{Proof of Proposition~\ref{P3.3}}
For simplicity, we will set $s = 0$.

First, let $X$ be a solution to the martingale problem related to $L$.
Since $h\in{\mathcal {D}}_L$ and $Lh=0$, we know that $Y=h(X)$ is an
$({\mathcal  F}_t)$-local martingale. In order to calculate its
bracket, we recall that $h^2\in{\mathcal {D}}_L$ and $Lh^2 = \sigma^2
(h')^ 2$ hold by Proposition~\ref{T212}(a). Thus,
\[
h^2(X_t) - \int_0^t (\sigma h')^2(X_s)\,ds
\]
is an $({\mathcal  F}_t)$-local martingale. This implies that
\[
[Y,Y]_t = \int_0^t (\sigma h')^2 (h^{-1}(Y_s))\,ds = \int_0^t
\tilde\sigma_h^2(Y_s)\,ds.
\]
Finally, $Y$ is a solution to the SDE (\ref{E3.3}) with respect to the
standard ${\mathcal {F}}_Y$-Brownian motion $W$ given by
\[
W_t = \int_0^t \frac{1}{\tilde\sigma_h(Y_s)} \, dY_s,
\]
where ${\mathcal  F}_Y$ is the canonical filtration generated by $Y$.

Now, let $Y=h(X)$ be a solution to (\ref{E3.3}) and let $f\in{\mathcal
{D}}_L$. Proposition~\ref{T212}(c) says that $\phi:= f\circ
h^{-1}\in{\mathcal {D}}_{L^0} \equiv C^2$, where
\begin{equation} \label{E3.4}
L^0 \phi= \frac{\tilde\sigma_h^2}{2} \phi'' = (Lf)\circ h^{-1}.
\end{equation}
We can therefore apply It\^{o}'s formula to evaluate $\phi(Y)$, which
coincides with $f(X)$. This gives
\[
\phi(Y_t) = \phi(Y_0) + \int_0^t \phi'(Y_s)\,dY_s + \tfrac
{1}{2}\int_0^t \phi''(Y_s)\,d[Y,Y]_s.
\]
Using $d[Y,Y]_s =\tilde\sigma_h^2(Y_s)\,ds$ and taking into
account (\ref{E3.4}), we conclude that
\begin{equation} \label{E3.5}
f(X_t) = f(X_0) + \int_0^t (f'\sigma)(X_s)\,dW_s + \int_0^t
Lf(X_s)\,ds.
\end{equation}
This establishes the proposition.
\end{pf*}

\begin{rem} \label{R3.5}
From Proposition~\ref{P3.3} in particular, we have the following.

Let $(\Omega, ({\mathcal  F}_t), P) $ be a filtered probability space
fulfilling the usual conditions. Let $x_0 \in{\mathbb R}$ and $X$ be a
solution to the martingale problem related to $L$ with initial
condition $x_0$. Then there exists a classical Brownian motion $(W_t)$
such that $X$ is a solution to the sharp martingale problem related to
$L$ with initial condition $x_0$.
\end{rem}

\begin{cor} \label{C3.5}
Let $X$ be a solution to the martingale problem related to $L$ with
initial condition $x_0$. Then map $\mathcal {A}$ admits a continuous
extension from ${\mathcal {D}}_L$ to $C^1$ with values in ${\mathcal
C}$ which we will again denote by $\mathcal {A}$. Moreover, ${\mathcal
{A}} (f)$ is a zero quadratic variation process for every $f\in C^1$.
\end{cor}

\begin{pf}
$\mathcal {A}$ has a continuous extension because of (\ref{E3.5}).
${\mathcal  A} (f)$ is a zero quadratic variation process because $X$
is a Dirichlet process with martingale part $\int_0^{\bolds{\cdot}}
\sigma(X_s)\,dW_s$ and because of Remark~\ref{R81}.
\end{pf}

\begin{rem} \label{R3.6}
The extension of (\ref{E3.5}) to $C^1$ gives
\begin{equation} \label{E3.6}
f(X_t)= f(X_0) + \int_0^t (f'\sigma)(X_s)\,dW_s + {\mathcal {A}}(f).
\end{equation}
Choosing $f=id$ in (\ref{E3.6}), we get
\[
X_t = X_0 + \int_0^t \sigma(X_s)\,dW_s + {\mathcal {A}}(id).
\]
\end{rem}

We will see that if there is a subspace $\nu_0$ of $C^0$ such that
Technical Assumption~\ref{A0} is verified, then the operator $A$ will
be extended to $\nu_0$. If $b$ is an element of that space, then it
will be possible to write $\hat L id = b$ and ${\mathcal {A}}(id) =
A(b) $. In that case, we will be able to indicate that $X$ is a
solution of the generalized SDE with diffusion coefficient $\sigma$ and
distributional drift $b'$.

A similar result to Proposition~\ref{P3.3} can be deduced for the case
of a transformation through function $k$ and the divergence-type
operator introduced at (\ref{Efk}).

\begin{prop} \label{PMardiv}
We consider the transformation $k$ and the PDE operator $L^1$ introduced
at \textup{(\ref{Efk})} and in Lemma \textup{\ref{T216}}, respectively.

A process $X$ solves the martingale problem related to $L$ with initial
condition $x_0$ at time $s$ if and only if $Z=k(X)$ solves the
martingale problem related to $L^1$ with initial condition $k(x_0)$ at
time~$s$.
\end{prop}

\begin{pf}
This is an easy consequence of Lemma \ref{T216}.
\end{pf}

Let $x_0 \in{\mathbb R}$, $y_0 = h(x_0)$.
Let $\sigma, b, \Sigma, h$ be as in Section~\ref{s4}.

We set $\tilde\sigma_h = (\sigma e^{-\Sigma}) \circ h^{-1}$.

From Proposition~\ref{P3.3}, we have the following.

\begin{cor} \label{CMP}
\textup{(i)} Strong existence (resp., pathwise uniqueness) holds for
$\mathit{SMP} (L, x_0)$ if and only if strong existence (resp., pathwise
uniqueness) holds for the SDE
\[
dY_t = \tilde\sigma_h (Y_r)\,dW_r
\]
with initial condition $Y_0 = h(x_0)$.

\textup{(ii)} An analogous equivalence holds for weak existence (resp.,
uniqueness in law).
\end{cor}

From Proposition~\ref{P3.3}, we can deduce two other corollaries
concerning the well-posedness of our martingale problem.

\begin{cor} \label{CWE} Under the same assumptions as the previous
corollary, $\mathit{MP}(L,x_0)$ admits weak existence and uniqueness in law.
\end{cor}

\begin{pf}
The statement follows from point (i) of Corollary \ref{CMP} and from
the fact that the SDE (\ref{E3.3}) admits weak existence and uniqueness
in law because \mbox{$\tilde\sigma_h > 0$;} see Theorem~5.7, Chapter~5
of~\cite{ks}, or \cite{esd}.
\end{pf}

\begin{rem} \label{RSem}
By Corollary 5.11 of \cite{frw2}, it is immediate to see that
the solution is a semimartingale for each initial condition if and only
if $\Sigma$ is locally of bounded variation.

If $L$ is in divergence form [see Remark~\ref{R2.6}(a) with
$\alpha= 1$], then the solution corresponds to the process constructed
and studied by, for instance, Stroock~\cite{str}.
\end{rem}

\begin{cor} \label{C3.3SS} Suppose that either
$(\sigma, b) \in(D^{{1}/{2}}, C^{{1}/{2}})$ or $(b, \sigma)
\in(D^{{1}/{2}}, C^{{1}/{2}})$ and, moreover, that
\textup{(\ref{Nonexplo})} is satisfied. Then $\mathit{MP}(L,x_0)$ admits strong
existence and pathwise uniqueness.
\end{cor}

\begin{pf}
In this case, $\Sigma$ is well defined [see Remark~\ref{R2.6}(d)] and
$\sigma$ belongs to $ D^{{1}/{2}}$. Since $h^{-1}$ is of class $C^1$,
$\tilde\sigma_h$ is H\"older continuous with parameter $\frac{1}{2}$.
The SDE (\ref{E3.3}) admits pathwise uniqueness because of Theorem
3.5(ii) of \cite{ry} and weak existence, again through Theorem~5.7 of
\cite{ks}. The Yamada--Watanabe theorem (see \cite{ks}, Corollary 3.23,
Chapter~5) also implies strong existence for~(\ref{E3.3}). The result
follows from point (i) of Corollary \ref{CMP}.
\end{pf}

\section{A significant stochastic differential equation with
distributional drift}\label{s6}

In this section, we will discuss the case where the martingale problem
is equivalent to a stochastic differential equation to be specified.
First, one would need to give a precise sense to the generalized drift
$\int_0^{\bolds{\cdot}} b'(X_s)\,ds$, $b$ being a continuous function.

We will introduce a property related to a general process $X$. First,
we consider the linear map $A^X\dvtx
\ell\rightarrow\int_0^{\bolds{\cdot}} \ell'(X_s)\,ds$ defined on
$C^1({\mathbb R})$ with values in~${\mathcal C}$.

\begin{defi} \label{Deltr}
Let $\nu_1$ be a topological F-space (or, eventually, an inductive
limit of $F$-spaces) which is a topological linear subspace of $C^0
({\mathbb R})$ and such that $\nu_1 \supset C^1 ({\mathbb R})$. We will
say that $X$ has \textit{extended local time regularity with
respect to $\nu_1$} if:
\begin{itemize}
\item $A^X$ admits a continuous extension to $\nu_1$, which will still
be denoted by the same symbol;

\item $\int_0^{\bolds{\cdot}} g(X)\,d^- A^X(\ell)$ exists for every $g
\in C^2$ and every $ \ell\in\nu_1.$
\end{itemize}
\end{defi}

\begin{rem} \label{RBY}
The terminology related to \textit{local time} is natural in this
context. To illustrate this, we consider a general continuous process
$X$ having a local time $(L_t(a), t \in[0,T], a \in{\mathbb R})$ with
respect to Lebesgue measure, that is, fulfilling the density occupation
identity
\[
\int_0^t \varphi(X_s)\,ds = \int_{\mathbb R}\varphi(a) L_t(a)\,da,
\qquad t \in[0,T],
\]
for every positive Borel function $\varphi$. $X$ trivially has
extended local time regularity, at least with respect to $\nu_1 = C^1
$.

Let $\ell\in C^1$. Suppose for a moment that $(L_t(a))$ is a
semimartingale in $a$, as is the case, for instance, if $X$ is a
classical Brownian motion. In that case, one would have
\[
\int_0^t \ell'(X_s)\,ds = \int_0^t \ell'(a) L_t(a)\,da = - \int
_{\mathbb R}\ell(a) L_t(da).
\]
Clearly, the rightmost integral can be extended continuously in
probability to any $\ell\in C^0$, which implies that $X$ also has
extended local time regularity related to $ \nu_1 = C^0$. We remark
that \cite{by} gives general conditions on semimartingales $X$ under
which $L_t(da)$ is a good integrator, even if $(L_t(a))$ is not
necessarily a semimartingale in $a$.
\end{rem}

\begin{defi} \label{DSol} Let $(\Omega, ({\mathcal  F}_t), P)$  a
filtered probability space, $(W_t)$ a classical $({\mathcal
F}_t)$-Brownian motion and $Z$ an ${\mathcal  F}_0$-measurable random
variable. A process $X$ will be called a $\nu_1$-\textit{solution} of
the SDE
\begin{eqnarray*}
dX_t &=& b'(X_t)\,dt + \sigma(X_t)\,dW_t,
\\
X_0 &=& Z,
\end{eqnarray*}
if:
\begin{itemize}
\item $X$ has the extended local time regularity with respect to
$\nu_1$;

\item $X_t = Z + \int_0^t \sigma(X_s)\,dW_s + A^X (b)_t$;

\item $X$ is a finite quadratic variation process.
\end{itemize}
\end{defi}

\begin{rem} \label{RSol}
Suppose that $b \in\nu_1$. If $\nu_1 \subset\nu_1'$, then a
$\nu _1'$-solution is also a $\nu_1$-solution.

The previous definition is also new in the classical case,
that is, when $b'$ is a continuous function. A $\nu_1$-solution with
$\nu_1 = C^1$ corresponds to a solution to the SDE in the classical
sense. On the other hand, a $\nu_1$-solution with $\nu_1$ strictly
including $C^1$ is a solution whose local time has a certain additional
regularity.
\end{rem}

Even in this generalized framework, it is possible to introduce the
notions of \textit{strong $\nu_1$-existence}, \textit{weak
$\nu_1$-existence}, \textit{pathwise $\nu_1$-uniqueness} and
\textit{$\nu_1$-uniqueness in law}. This can be done similarly as in
Definition \ref{D115} according to whether or not the filtered
probability space with the classical Brownian motion is fixed a priori.

\begin{lemma} \label{L3.9}
We suppose that Technical Assumption \textup{\ref{A0}} is satisfied. If
$X$ is a solution to a martingale problem related to a PDE operator
$L$, then it has extended local time regularity with respect to $\nu_1
= \nu_0$.
\end{lemma}

\begin{pf}
Let $ \ell\in C^1$. Since $X$ solves the martingale problem with
respect to $L$, setting $f = T\ell$, it follows that
\begin{eqnarray*}
A^X (\ell)_t &=& \int_0^t \ell'(X_s)\,ds = \int_0^t Lf (X_s)\,ds
\\
&=& f(X_t) - f(X_0) - \int_0^t f'(X_s) \sigma(X_s)\,dW_s.
\end{eqnarray*}
Continuity of $T$ on $\nu_0$ implies that $A^X$ can be extended to
$\nu_0$.

Now, let $ \ell\in\nu_0$ and $f = T\ell\in C^1$. Since $f(X)$ equals a
local martingale plus $A^X(\ell)$, it remains to show that
\begin{equation} \label{E3.8}
\int_0^{\bolds{\cdot}} g(X)\,d^- f(X)
\end{equation}
exists for any $g \in C^2$. Integrating by parts, the previous integral
(\ref{E3.8}) equals
\[
(g f)(X_{\bolds{\cdot}}) - (gf)(X_0) - \int_0^{\bolds{\cdot}} f(X)\,d^-
g(X) - [f(X),g(X)].
\]
Remark~\ref{R1.2}(b), (f) shows that the rightmost term member is
well defined.
\end{pf}

\begin{lemma} \label{L3.10}
Let $X$ be a process having extended local time regularity with
respect to some $F$-space (or inductive limit) $\nu_1$. Suppose that
for fixed $g \in C^1$, the application $\ell\rightarrow g \ell$ is
continuous from $\nu_1$ to $\nu_1$. Then for every $g \in C^2$ and
every $ \ell\in\nu_1$, we have
\begin{equation} \label{E3.9}
\int_0^{\bolds{\cdot}} g(X)\,d^- A^X(\ell) = A^X(\Phi(g,\ell)),
\end{equation}
where
\begin{equation} \label{E3.10}
\Phi(g,\ell)(x) = (g\ell)(x) - (g\ell)(0) -\int_0^x (\ell g')(y)\,dy.
\end{equation}
\end{lemma}

\begin{pf}
The Banach--Steinhaus-type Theorem \ref{TBS} implies that for every
\mbox{$g\in C^2$},
\begin{equation} \label{E3.11}
\ell\mapsto\int_0^{\bolds{\cdot}} g(X)\,d^-A^X(\ell)
\end{equation}
is continuous from $\nu_1$ to $\mathcal {C}$. In fact, expression
(\ref{E3.11}) is the u.c.p. limit of
\[
\lim_{\varepsilon\to0+} \int_0^{\bolds{\cdot}} g(X_s)
\frac{A^X(\ell)_{s+\varepsilon}-A^X(\ell)_s}{\varepsilon}\,   ds.
\]
Note that $\Phi$ is a continuous bilinear map from $C^1 \times\nu_1$ to
$\nu_1$. Since $A^X\dvtx  \nu_1 \rightarrow{\mathcal  C}$ is
continuous, the mapping $\ell\rightarrow A^X(\Phi(g,\ell))$ is also
continuous from $\nu_1$ to ${\mathcal  C}$. In order to conclude the
proof, we need to check identity (\ref{E3.9}) for $\ell\in C^1$. In
that case, since
\[
\Phi(g,\ell)(x) = \int_0^x (g \ell')(y)\,dy,
\]
both sides of (\ref{E3.9}) equal
\[
\int_0^{\bolds{\cdot}}(g\ell')(X_s)\,ds.
\]\upqed
\end{pf}

We will now explore the relation between the martingale problem
associated with $L$ and the stochastic differential equations with
distributional drift.

\begin{prop} \label{P3.11} Let $x_0 \in{\mathbb R}$.
Suppose that $L$ fulfills Technical Assumption~\textup{\ref{A0}}. Let
$(\Omega, ({\mathcal  F}_t), P)$ be a filtered probability space
fulfilling the usual conditions and let $(W_t)$ be a classical
$({\mathcal F}_t)$-Brownian motion.

If $X$ solves the sharp martingale problem with respect to $L$ with
initial condition~$x_0$, then $X$ is a $\nu_0$-solution to the
stochastic differential equation
\begin{eqnarray} \label{Egen}
dX_t &=& b'(X_t)\,dt + \sigma(X_t)\,dW_t,
\nonumber
\\[-8pt]
\\[-8pt]
\nonumber X_0 &=& x_0.
\end{eqnarray}
\end{prop}

\begin{rem}\label{R2.21}
In particular, if $L$ is close to divergence type, as in
Example~\ref{4ex}(ii), then $X$ is a $C^0$-solution to the previous
equation with $b=\frac{\sigma^2}{2} + \beta- \frac{\sigma^2(0)}{2}$.
\end{rem}

\begin{pf}
Let $X$ be a solution to the martingale problem related to $L$. We
know, by Lemma \ref{L3.9}, that $X$ has extended local time regularity
with respect to $\nu_1$. On the other hand, by Remark~\ref{R3.4}, $X$
is a finite quadratic variation process. It remains to show that
\begin{equation} \label{E3.12}
X_t = X_0 + \int_0^t \sigma(X_s)\,dW_s + A^X(b)_t.
\end{equation}
Let $\ell\in C^1$ and set $f = T^1 \ell$. By definition of a sharp
martingale problem, we have

\begin{equation} \label{E3.13}
T^1 \ell(X_t) = T^1 \ell(X_0) + \int_0^t ((T^1\ell)' \sigma)
(X_s)\,dW_s + A^X (\ell)_t.
\end{equation}
According to Remark~\ref{R2.21a}(i) concerning the continuity of the
map $T^1\dvtx  \nu_0 \rightarrow C^1$, previous expression can be
extended to any $\ell\in\nu_0$.

By Remark~\ref{R2.21a}(ii), $\ell= b\in\nu_0$ and $f = T^1 \ell= id$.
Replacing this in (\ref{E3.13}), we obtain
\[
X_t = x_0 + \int_0^t \sigma(X_s)\,dW_s + A^X (b).
\]
Since $X_0 = Z$, the proof is complete.
\end{pf}

\begin{cor} \label{C3.11}
Let $x_0 \in{\mathbb R}$. Suppose that $L$ fulfills Technical
Assumption~\textup{\ref{A0}}. If $\mathit{MP}(L,x_0)$ [resp. $\mathit{SMP}(L,x_0)$] admits weak (resp.,
strong) existence, then the SDE \textup{(\ref{Egen})} also admits weak
(resp., strong) existence.
\end{cor}

\begin{pf}
The statement concerning strong solutions is obvious. Concerning weak
solutions, let us admit the existence of a filtered probability space,
where there is a solution to the martingale problem with respect to $L$
with initial condition $x_0$. Then according to Remark~\ref{R3.5}, this
solution is also a solution to a sharp martingale problem and the
result follows.
\end{pf}

If $X$ is some $\nu_1$-solution to (\ref{E3.12}), is it a solution to
the (sharp) martingale problem related to some operator $L$? This is a delicate
question.
In the following proposition, we only provide the converse
of Proposition~\ref{P3.11} as a partial answer.

\begin{prop} \label{C3.12}
Suppose that the PDE operator $L$ fulfills Technical Assumption
\textup{\ref{A0}}. Let $(\Omega, ({\mathcal  F}_t), P)$ be a filtered
probability space fulfilling the usual conditions and let $(W_t)$ be a
classical $({\mathcal F}_t)$-Brownian motion. Let $X$ be a
progressively measurable process.

$X$ solves the sharp martingale problem related to $L$ with respect to
some initial condition $x_0$ if and only if it is a $\nu_0$-solution to
the stochastic differential equation
\begin{eqnarray} \label{E3.15}
dX_t &=& b'(X_t)\,dt + \sigma(X_t)\,dW_t, \nonumber
\\[-8pt]
\\[-8pt]
\nonumber  X_0 &=& x_0 .
\end{eqnarray}
\end{prop}

\begin{cor} \label{CSWEqu} Let $x_0 \in{\mathbb R}$.
Suppose that $L$ fulfills Technical Assumption \textup{\ref{A0}}. Then
weak existence and uniqueness in law (resp., strong existence and
pathwise uniqueness) hold for equation \textup{(\ref{E3.15})} if and
only if the same holds for~$\mathit{MP}(L,x_0)$ [resp. $\mathit{SMP}(L,x_0)$].
\end{cor}

\begin{pf*}{Proof of Proposition~\ref{C3.12}}
Suppose that $X$ is a $\nu_0$-solution to~(\ref{E3.15}). Then it is a
finite quadratic variation process. Let $f\in C^3$. Since $X$
solves~(\ref{E3.12}) and $\int_0^{\bolds{\cdot}} f'(X_s)\,d^-X_s$ always
exists by the classical It\^{o} formula [see Remark~\ref{R1.2}(e) of
Section~\ref{s1}], we know that $\int_0^{\bolds{\cdot}} f'(X)\,d^-
A^X(b)$ also exists and is equal to $\int_0^{\bolds{\cdot}} f'(X)\,d^-X
- \int_0^{\bolds{\cdot}} (f'\sigma)(X)\,dW.$ Therefore, this It\^{o}
formula says that
\begin{eqnarray*}
f(X_t) & = & f(X_0) + \int_0^t f'(X_s) \sigma(X_s)\,dW_s
+ \int_0^t f'(X)\,d^- A^X(b)
\\
&& {}+\tfrac{1}{2} \int_0^t f''(X_s)\sigma^2(X_s)\,ds
\end{eqnarray*}
holds.

By Lemma \ref{L3.10}, the linearity of mapping $A^X$ and (\ref{E213a}),
we obtain
\begin{eqnarray*}
&& \int_0^{t} f'(X)\,d^- A^X(b) + \frac{1}{2}\int_0^t (f'' \sigma
^2)(X_s)\,ds
\\
&&\qquad =  A^X (\Phi(f',b) )_t + \frac{1}{2}\int_0^t (f'' \sigma
^2)(X_s)\,ds
\\
&&\qquad =  \int_0^t \biggl(\frac{\sigma^2}{2} - b\biggr)(X_s)
f''(X_s)\,ds + A^X(b f') = A^X(\hat L f).
\end{eqnarray*}

This shows that
\begin{equation} \label{E3.14}
f(X_t) - f(X_0) - \int_0^t (f' \sigma)(X_s)\,dW_s = A^X(\hat L f)
\end{equation}
for every $f \in C^3$. In reality, it is possible to show the previous
equality for any $f \in C^2$. In fact, the left-hand side extends
continuously to $C^2$ and even to $C^1$. The right-hand side is also
allowed to be extended to $C^2$ for the following reason. For $f\in
C^2$, let $(f_n)$ be a sequence of functions in $C^3$ converging to $f$
when $n \rightarrow\infty$, according to the $C^2$ topology. In
particular, the convergence also holds in $C^1_{\nu_0}$. Since $\hat L$
is continuous with respect to the $C^1_{\nu_0}$ topology with values in
$\nu_0$, we have $\hat L f_n \rightarrow\hat L f$ in $\nu_0$. Finally,
$ A^X(\hat L f_n) \rightarrow A^X(\hat L f)$ u.c.p. because of the
extended local time regularity with respect to $\nu_0$.

We will, in fact, use the validity of (\ref{E3.14}) for $f \in C^2$
with $f(0) = 0$ and  $x_1 = f'(0)$ and $\ell= \hat L f.$ According
to Technical Assumption \ref{A0}(iv), we have $f = T^{x_1} \ell$.
Therefore, (\ref{E3.14}) gives
\[
T^{x_1}\ell(X_t) = T^{x_1} \ell(X_0) + \int_0^t ((T^{x_1} \ell)'
\sigma)(X_s)\,dW_s + A^X (\ell).
\]
Again using extended local time regularity with respect to $\nu _0$ and
the continuity of $T^{x_1}$, we can state the validity of the previous
expression for each $\ell\in\nu_0$ with $\ell(0)=0$, in particular, for
$\ell\in C^1$ with $\ell(0)=0$. But in this case, for any $f
\in{\mathcal  D}_L$ with $f(0) = 0$ and $ \ell ' = Lf$, we obtain
\[
f(X_t) = f(X_0) + \int_0^t (f' \sigma)(X_s)\,dW_s + \int_0^t
Lf(X_s)\,ds.
\]
This shows the validity of the identity in Definition \ref{DSmp} for $f
\in{\mathcal  D}_L$ and that $f(0) = x_0$ and $x_0 = 0$. If $x_0 \neq0$, we
replace $f$ by $f - x_0$ in the previous identity and use the fact that
$L (f - x_0) = L f$ for any $ f \in{\mathcal  D}_L$.

It follows that $X$ fulfills a sharp martingale problem with respect
to~$L$.

This shows the reversed sense of the statement. The direct implication
was proven in Proposition~\ref{P3.11}.
\end{pf*}

\begin{cor} \label{CEqStrong}
We suppose that $\sigma\in D^{{1}/{2}}$ and $b \in C^{{1}/{2}}$, or
$\sigma\in C^{{1}/{2}}$ and $b \in D^{{1}/{2}}$, with
conditions~\textup{(\ref{Nonexplo})}. We set $\nu_0 = D^{{1}/{2}}$.

Then equation (\ref{E3.15}) admits $\nu_0$-strong existence and
pathwise uniqueness.
\end{cor}

\begin{pf}
The result follows from Corollaries \ref{CSWEqu}
and \ref{C3.3SS}.
\end{pf}

\section{About $C^0_b$-generalized solutions of parabolic
equations}\label{s7}

In this section, we want to discuss the related parabolic Cauchy
problem with final condition, which is associated with our stochastic
differential equations with distributional drift.

We will adopt the same assumptions and conventions as in
Section~\ref{s4}. We consider the formal operator ${\mathcal  L}=
\partial_t + L$, where $L$ will hereafter act on the second
variable.

\begin{defi} \label{D21}
Let $\lambda$ be an element of $C^0_b( [0,T] \times{\mathbb R}) $ and
let $u^0 \in C^0_b({\mathbb R})$. A function $u\in C^0_b([0,T]
\times{\mathbb R})$ will be said to be a \textit{$C^0_b$-generalized
solution} to
\begin{eqnarray} \label{F25}
{\mathcal  L}u &=& \lambda,
\nonumber
\\[-8pt]
\\[-8pt]
\nonumber u(T,\cdot) &=& u^0,
\end{eqnarray}
if the following are satisfied:
\begin{longlist}[(iii)]
\item[(i)] for any sequence $(\lambda_n)$ in $C^{0}_b([0,T]
\times{\mathbb R})$ converging to $\lambda$ in a bounded
way,\looseness=1

\item[(ii)] for any sequence $(u_n^0)$ in $C_b^0({\mathbb R})$
converging in a bounded way to $u^0$,

\item[(iii)] such there are classical solutions $(u_n)$ in $C^{0}_b
([0,T]\times{\mathbb R})$ of class\break $C^{1,2} ([0,T[\times{\mathbb
R})$ to ${\mathcal  L}_n u_n = \lambda_n$, $u_n(T,\cdot) = u^0_n$,
\end{longlist}
then $(u_n)$ converges in a bounded way to $u$.
\end{defi}

\begin{rem} \label{R21}
(a) $u$ is said to solve ${\mathcal  L}u =\lambda$ if there exists $u^0
\in C^0_b({\mathbb R})$ such that (\ref{F25}) holds.
{\smallskipamount=0pt
\begin{longlist}[(a)]
\item[(b)] The previous definition depends in principle on the
mollifier, but it could be easily adapted so as not to depend on it.

\item[(c)] The regularized problem admits a solution: if $u_n^0 \in
C_b^3({\mathbb R})$ and $ \lambda_n \in C^{0,1}_b([0,T] \times{\mathbb
R})$, then there is a classical solution $u_n$ in
$C^{1,2}([0,T]\times{\mathbb R})$ of
\begin{eqnarray*}
{\mathcal  L}_n v &=& \lambda_n,
\\
v(T,\cdot) &=& u^0_n.
\end{eqnarray*}
For this, it suffices to apply Theorem \textup{5.19}
of~\textup{\cite{l}}.
\end{longlist}}
\end{rem}

We now state a result concerning the case when the operator $L$ is
classical. Even if the next proposition could be stated when the drift
$b'$ is a continuous function, we will suppose it to be zero. In fact,
it will later be applied to $L = L^0$.

\begin{prop} \label{P39}
We suppose that $b = 0$. Let $\varphi, \varphi_n \in C^0_b({\mathbb
R}), g, g_n \in C_b^0 ([0,T]\times{\mathbb R}), n \in{\mathbb N},$ such
that $\varphi_n \longrightarrow\varphi$, $ g_n \longrightarrow g$ in a
bounded way on ${\mathbb R}$ and $[0,T]\times{\mathbb R}$.

Let $\sigma$ be a strictly positive real continuous function.

Suppose that there exist $u_n\in C^{1,2}([0,T[\,\times\,{\mathbb R})
\cap C_b^{0}([0,T]\times {\mathbb R}) $ such that
\begin{eqnarray*}
{\mathcal  L}_n u_n &=& g_n,
\\
u_n (T,\cdot) &=&\varphi_n.
\end{eqnarray*}
Then $(u_n)$ will converge to $u \in C^0_b([0,T]\times{\mathbb R})$ in
a bounded way, where the function $u$ is defined by
\begin{equation} \label{F311ter}
u (s,x) = {\mathbb E} \biggl( \varphi(Y_T^{s,x}) + \int_s^T
g(r,Y_T^{r,x})\,dr
 \biggr),
\end{equation}
where $Y= Y^{s,x}$ is the unique solution (in law) to
\begin{equation} \label{F310quater}
Y_t = x + \int_s^t \sigma(X_r)\,dW_r
\end{equation}
and where $(W_t)$ is a classical Brownian motion on some suitable
filtered probability space.
\end{prop}

\begin{rem} \label{R310}
Usual It\^o calculus implies that
\begin{equation} \label{F311}
u_n(s,x) = {\mathbb E} \biggl( \varphi_n (Y_T^{s,x}(n)) + \int_s^T
g_n(r,Y_T^{r,x}(n))\,dr  \biggr),
\end{equation}
where $Y(n) = Y^{s,x}(n)$ is the unique solution in law to the problem
\begin{equation} \label{F310quinque}
Y_t(n) = x +\int_s^t \sigma_n (Y_r(n))\,dW_r.
\end{equation}

Theorem \textup{5.4} (Chapter~\textup{5} of \textup{\cite{ks}})
affirms that it is possible to construct a solution (unique in law) $Y
= Y^{s,x}$ to the SDE (\ref{F310quater}) \textup{[}resp., $Y(n) =
Y^{s,x} (n) $ to \textup{(\ref{F310quinque})]}.\vadjust{\goodbreak}

Suppose that $L$ is a classical PDE operator. Let $u \in
C^{1,2}([0,T[\times{\mathbb R})$ be bounded and continuous on
$[0,T]\times{\mathbb R}$. Again, It\^o calculus shows that $u$ can be
represented by (\ref{F311ter}) and (\ref{F310quater}). In particular, a
classical solution $u$ to ${\mathcal  L}u = g$ is also a
\mbox{$C^0_b$-}generalized solution.
\end{rem}

\begin{pf*}{Proof of Proposition~\ref{P39}}
We fix $s \in[0,T]$, $x \in{\mathbb R}$. Using the Engelbert--Schmidt
construction (see, e.g., the proof of Theorem 5.4, Chapter~5 and 5.7
of~\cite{ks}), it is possible to construct a solution $Y = Y^{s,x}$ of
the SDE on some fixed probability space which solves (\ref{F310quater})
with respect to some classical Wiener process~$(W_t)$. We set $s = 0$
for simplicity. The procedure is as follows. We fix a standard Brownian
motion $(B_t)$ on some fixed probability space  one set
\[
R_t := \int_0^t \frac{du}{\sigma^2(x + B_u)}.
\]
$R$ is a.s. a homeomorphism on ${\mathbb R}_+$ and we define $A$ as the
inverse of $R$. A solution $Y$ will be then given by $Y_t = x +
B_{A_t}$; in fact, it is possible to show that the quadratic variation
of the local martingale $Y$ is
\[
\langle Y,Y\rangle _t = \int_0^t \sigma^2(Y_s)\,ds.
\]
The Brownian motion $W$ is constructed a posteriori and  is adapted to
the natural filtration of $Y$ by setting $ W_t = \int_0^t
\frac{dY_s}{\sigma(Y_s)}.$

So, on the same probability space, we can set $Y_t(n) = x +
B_{A_t(n)}$, $A(n)$ being the inverse of $R(n)$, where $R(n)_t :=
\int_0^t \frac{du}{\sigma_n^2(x + B_u)}$.

Consequently, on the same probability space, we construct $Y_t(n) = x +
B_{A_t(n)} $, where $A(n)$ is the inverse of $R(n)$ and $R(n)_t :=
\int_0^t \frac{du}{\sigma_n^2(x + B_u)}$. $Y(n)$ solves
equation~(\ref{F310quinque}) with respect to a Brownian motion
depending on~$n$.

By construction, the family $Y_T^{s,x}(n)$ converges a.s. to
$Y_T^{s,x}$. Using Lebesgue dominated convergence theorems and the
bounded convergence of $(\varphi_n)$ and $(g_n)$, we can take the limit
when $n \rightarrow\infty$ in expression (\ref {F311}) and obtain the
desired result.
\end{pf*}

\begin{rem} \label{R311a}
In particular, the corresponding laws of random variables
$(Y^{s,x}(n))$ are tight.
\end{rem}

Again, we will adopt the same conventions as in Section~\ref{s4}.

We set $\sigma_h =\sigma h'$.  $L^0$ is the classical operator defined
at (\ref{E2.12}). Let us consider ${\mathcal  L}^0 = \partial_t + L^0$
as a formal operator.

\begin{cor} \label{C311}
Let $g \in C^0_b([0,T] \times{\mathbb R})$, $\varphi\in C^0_b({\mathbb
R})$. There is a \mbox{$C_b^0$-}gene\-ralized solution $u$ to
${\mathcal L}^0 u=g $, $ u(T,\cdot) =\varphi$. This solution is unique
and is given by~\textup{(\ref{F311ter}).}
\end{cor}

We now return to the original PDE operator ${\mathcal  L}$ with
distributional drift. We again denote by $h$ the same application defined in
Section~\ref{s5} and discuss existence and uniqueness of
$C^0_b$-generalized solutions of related parabolic Cauchy problems.

A useful consequence of Proposition~\ref{P39} is the following.
\begin{theo} \label{T313}
For $\varphi\in C^0([0,T] \times{\mathbb R})$ or $C^0({\mathbb R})$, we
again set $\tilde{\varphi} = \varphi\circ h^{-1}$ according to the
conventions of Section~\textup{\ref{s2}}. Again, we consider ${\mathcal
L}^0 = \partial_t + L^0$ as a formal operator.

Let $\lambda\in C^0_b([0,T] \times{\mathbb R})$, $ u^0 \in
C^0_b({\mathbb R}).$

There is a unique solution $u \in C^0_b([0,T] \times{\mathbb R})$ to
\begin{eqnarray} \label{F316}
{\mathcal  L}u &=& \lambda, \nonumber \\[-8pt]
\\[-8pt]
\nonumber u(T,\cdot) &=& u^0.
\end{eqnarray}
Moreover, $\tilde{u}$ solves
\begin{eqnarray} \label{F316a}
{{\mathcal  L}}^0 \tilde{u} &=& \tilde\lambda,
\nonumber
\\[-8pt]
\\[-8pt]
\nonumber \tilde{u} (T,\cdot) &=& \tilde{u}^0.
\end{eqnarray}
\end{theo}

\begin{pf}
In accordance with Section~\ref{s4}, let $(h_n)_{n\in{\mathbb N}}$ be
an approximating sequence which is related to $L h = 0$. Let us
consider the PDE operators ${\mathcal  L}_n$ defined at (\ref{E2.2}).
Let $(\lambda_n)_{n\in{\mathbb N}}$ be a sequence in
$C^0_b([0,T)\times{\mathbb R})$ such that $\lambda_n\to\lambda$, $u_n^0
\to u^0$ in a bounded way and for which there are classical solutions
$u_n$ of
\begin{eqnarray*}
{\mathcal  L}_n u_n &=& \lambda_n,
\\
u_n(T,\cdot) &=& u^0_n.
\end{eqnarray*}
We recall that those sequences always exist because of
Remark~\ref{R21}(c).

We set
\[
g_n = \lambda_n \circ h_n^{-1},\qquad  \varphi_n = \varphi\circ
h_n^{-1},\qquad   v_n = u_n \circ h_n^{-1}.
\]
By Lemma $\ref{L28}$, we have
\begin{eqnarray*}
{{\mathcal  L}}_n^0 v_n &=& g_n,
\\
v_n(T,\cdot) &=& \varphi_n,
\end{eqnarray*}
where
\[
{{\mathcal  L}}_n^0 \varphi(t,y)=\partial_t \varphi(t,y)+
{\sigma}^2_{h_n} \circ h_n^{-1} (t,y)\, \partial_{\mathit{xx}}^2
\varphi(t,y).
\]
By Proposition~\ref{P39}, and Corollary \ref{C311}, $v_n \to\tilde{u}$
in a bounded way, where
\begin{eqnarray*}
{{\mathcal  L}}^0 \tilde{u} &=& \tilde{\lambda},
\\
\tilde{u}(T,\cdot) &=& \tilde{u}_0.
\end{eqnarray*}
This concludes the proof of the proposition.
\end{pf}

We now discuss how $C^0_b$-generalized solutions are transformed under
the action of the function $k$ introduced at (\ref{Efk}). A similar
result to Lemma \ref{T216} for the elliptic case is the following.

\begin{prop} \label{P313}
For $\varphi\in C^0([0,T] \times{\mathbb R})$ or $C^0({\mathbb R})$, we
set $\bar{\varphi} = \varphi\circ k^{-1}$. We set $\sigma_k =\sigma k'$
and consider the formal operator
\[
{{\mathcal  L}}^1 f = \partial_t f +\tfrac{1}{2}\bar{\sigma}_k^2\,
\partial ^2_{\mathit{xx}} f + \tfrac{1}{2} (\bar{\sigma}_k^2)'\, \partial_x f.
\]
Informally, we can write
\[
{{\mathcal  L}}^1 f = \partial_t f + \tfrac{1}{2}\partial_x
(\bar{\sigma }_k^2 \,\partial_x f).
\]

Let $\lambda\in C^0_b([0,T] \times{\mathbb R})$, $ u^0 \in
C^0_b({\mathbb R}).$

Let $u$ be the unique $C^0_b$-generalized solution in $ C^0_b([0,T]
\times{\mathbb R})$ to
\begin{eqnarray} \label{F316bis}
{\mathcal  L}u &=& \lambda,
\nonumber
\\[-8pt]
\\[-8pt]
\nonumber u(T,\cdot) &=& u^0.
\end{eqnarray}
Then $\bar{u}$ solves
\begin{eqnarray*}
{{\mathcal  L}}^1 \bar{u} &=& \bar\lambda,
\\
\bar{u} (T,\cdot) &=& \bar u^0.
\end{eqnarray*}

\end{prop}
\begin{pf}
Let $v$ be the unique solution to
\begin{eqnarray*}
{{\mathcal  L}}^1 v &=& \bar\lambda,
\\
v (T,\cdot) &=& \bar u^0,
\end{eqnarray*}
which exists because of Theorem \ref{T313}, taking ${\mathcal  L}=
{\mathcal L}^1$.

We define $H\dvtx  {\mathbb R}\rightarrow{\mathbb R}$ such that
\[
H(0) = 0,\qquad   H'(z) = \frac{1}{\sigma_k^2} (z).
\]
Again, (\ref{Nonexplo}) implies that $H$ is bijective on ${\mathbb R}$.
This case corresponds to example~(a) in Remark~\ref{R2.6} with $\alpha
= 1$.

We set $\tilde v = v \circ H^{-1}$. Again, by Theorem \ref{T313}, we
have
\begin{eqnarray*}
{{\mathcal  L}}^{0,1} \tilde v &=& \bar\lambda\circ H^{-1},
\\
\tilde v (T,\cdot) &=& u^0 \circ(k^{-1} \circ H^{-1}),
\end{eqnarray*}
where ${{\mathcal  L}}^{0,1} f = \frac{a^2}{2}\,
\partial_{\mathit{xx}}^2 f$ and
\[
a = (\sigma_k H') \circ H^{-1} = \frac{1}{\sigma_k} \circ H^{-1}.
\]
Since
\[
\sigma_k = (\sigma k') \circ k^{-1} = \frac{e^{\Sigma}}{\sigma}
\circ k^{-1},
\]
this yields
\[
a = (\sigma e^{-\Sigma}) \circ(H \circ k)^{-1}.
\]
On the other hand, $H \circ k = h$ since
\begin{eqnarray*}
H \circ k (0) & = & 0 = h(0),
\\
(H \circ k (x))' & = & H'(k (x)) k'(x) = \frac{1}{\sigma_k^2} k'(x) =
\frac{1}{\sigma^2 k'} = e^{-\Sigma}  = h'.
\end{eqnarray*}
We can therefore conclude that ${{\mathcal  L}}^{0,1} \equiv{\mathcal
L}^0$. Since problem (\ref{F316a}) has a unique solution, $\tilde v =
\tilde u$, where $u$ solves (\ref{F316}) and $\tilde u = u \circ
h^{-1}. $ Finally,
\[
v = \tilde v \circ H = \tilde u \circ H = u \circ H \circ h^{-1}  = u
\circ k^{-1} = \bar u.
\]\upqed
\end{pf}

\begin{prop} \label{Probrep}
The unique $C^0_b$-generalized solution to \textup{(\ref{F316})} admits
a probabilistic representation in the sense that
\begin{equation} \label{F311quinque}
u (s,x) = {\mathbb E} \biggl( u^0 (X_T^{s,x}) + \int_s^T \lambda
(r,X_T^{r,x})\,dr  \biggr),
\end{equation}
where $X^{s,x}$ is the solution to the martingale problem related to
$L$ at time $s$ and point $x$.
\end{prop}

\begin{pf}
The result follows from Theorem \ref{T313}, Corollary \ref{C311} and
Proposition~\ref{P3.3}, which collectively imply the following. If $X$
is a solution to the martingale problem related to $L$ at point $x$ at
time $s$, then $Y = h(X)$ solves the stochastic differential equation
(\ref{E3.3}) with initial condition $h(x)$ at time $s$.
\end{pf}

\section{Density of the associated semigroups}\label{s8}

We now discuss the existence of a density law for the solutions
$X^{s,x}$ of the martingale problem related to $L$. First,  we suppose
that $L$ is an operator in divergence form with $Lf =
(\frac{\sigma^2}{2} f')' $ and that there are positive constants such
that $c \le\sigma^2 \le C$. We will say, in this case, that $L$
\textit{has the Aronson form}. This terminology refers to the
fundamental paper \cite{ar} concerning exponential estimates of
fundamental solutions of nondegenerate parabolic equations. We begin
with some properties (partly classical) stated in \cite{frw2}. We
observe that point (ix) is slightly modified with respect to
\cite{frw2}, but this new configuration can be immediately deduced from the
proof in \cite{frw2}. This preparatory work will be applied to the
operator $L^1$ introduced in (\ref{E2.21}).

\begin{lemma} \label{L5.1}
We suppose that $0 < c\le\sigma^2 \le C$. Let $\sigma_n$, $n\in\mathbb
{N}$, be smooth functions such that $0 < c\le\sigma_n^2\le C$ and
$\sigma_n^2\to \sigma^2$ in $C^0$, as at the beginning of
Section~\ref{s4}. We set $L_n g =(\frac{\sigma _n^2}{2}g')'$. There
exists a family of probability measures $(\nu_t(dx,y),
t\ge0,y\in\mathbb{R})$ [resp., $(\nu^n_t(dx,y), t\ge0,y\in\mathbb{R})$]
enjoying the following properties:

\begin{longlist}[(viii)]
\item[(i)] $\nu_t(dx,y) = p_t(x,y)\,dx$, $\nu_t^n(dx,y) =
p_t^n(x,y)\,dy$;

\item[(ii)] (Aronson estimates) there exists $M>0$, depending only on
constants~$c$, $C$, with
\[
\frac{1}{M\sqrt{t}} \exp \biggl(-\frac{M|x-y|^2}{t} \biggr) \le
p_t(x,y) \le\frac{M}{\sqrt{t}} \exp \biggl(-\frac{|x-y|^2}{Mt} \biggr);
\]

\item[(iii)] we have
\begin{equation} \label{E5.2}
\partial_t \nu_t(\cdot,y) = L\nu_t(\cdot,y),\qquad
\nu_0(\cdot,y) = \delta_y
\end{equation}
and
\[
\partial_t \nu^n_t(\cdot,y) = L_n\nu^n_t(\cdot,y),\qquad
\nu^n_0(\cdot,y) = \delta_y,
\]
where $\nu$ (resp., $\nu^n$) is called the \textup{fundamental
solution} related to the previous parabolic linear equation;

\item[(iv)] we have
\begin{eqnarray*}
\partial_t \nu_t (x,\cdot) &=& L \nu(x,\cdot),
\\
\partial_t \nu_t^n (x,\cdot) &=& L_n \nu^n (x, \cdot);
\end{eqnarray*}

\item[(v)] the map $(t,x,y)\mapsto p_t(x,y)$ is continuous from
$]0,\infty [\times \mathbb{R}^2$ to $\mathbb{R}$;

\item[(vi)] the $p^n$ are smooth on $]0,\infty[ \times\mathbb{R}^2$;

\item[(vii)] we have $\lim_{n\to\infty} p_t^n(x,y) = p_t(x,y)$
uniformly on each compact subset of $]0,\infty[\times\mathbb{R}^2$;

\item[(viii)] $p_t(x,y)= p_t(y,x)$ holds for every $t>0$ and every
$x,y\in \mathbb{R}$;

\item[(ix)] $ \int_0^T \sup_{y}  ( \int_{\mathbb{R}} |\partial _x p_t
(x,y) |^2\, dx  )^{{1}/{2}} \,dt < \infty.$
\end{longlist}
\end{lemma}

The previous lemma allows us to establish the following.

\begin{theo} \label{TDistr} Let $Z^{s,x}$ be the solution to the
martingale problem related to $L$ at time $s$ and point $x$. Suppose
that $L$ to be of divergence type, having the Aronson form. Then there
is fundamental solution $\nu_t = r_t(x,y)$ of
\[
\partial_t \nu_t (\cdot,y) = L\nu_t(\cdot,y),\qquad
\nu_0(\cdot,y) = \delta_y,
\]
with the following properties:
\begin{longlist}[(ii)]
\item[(i)] letting $g \in C^0_b([0,T] \times{\mathbb R})$, $\varphi\in
C_b({\mathbb R})$, the $C^0_b$-generalized solution $u$ to ${\mathcal
L}u=g$, $ u(T,\cdot) =\varphi$, is given by
\begin{equation} \label{F312terd}
\qquad u(s,x) = \int_{\mathbb R}\varphi(y) r_{T-s}(x,y)\,dy + \int_s^T\, dr
\int_{\mathbb R}g(r,y) r_{T-r} (x,y)\,dy;
\end{equation}

\item[(ii)] the law of $Z_T^{s,x}$ has $r_{T-s}(x, \cdot)$ as density
with respect to Lebesgue measure.
\end{longlist}
\end{theo}

\begin{pf}
Let $(r_t^n (x,y))$ be the fundamental solution corresponding to the
parabolic equation associated with $L_n f (x) = (\frac{\sigma_n^2
f'}{2})' $, as introduced in Section~\ref{s4}. We observe that
$(\sigma_n^2) $ converges in a bounded way to $ \sigma^2$.
\begin{longlist}[(ii)]
\item[(i)] We define
\begin{equation} \label{F312sex}
\qquad u_n(s,x) = \int_{\mathbb R}\varphi(y) r^n_{T-s}(x,y)\,dy + \int_s^T\,
dr \int_{\mathbb R}g(r,y) r^n_{T-r} (x,y)\,dy .
\end{equation}
Points (vi) and (ii) of Lemma \ref{L5.1} imply that functions $u_n$
belong to\break $C^{1,2}( [0,T[ \times{\mathbb R})$, so they are
classical solutions to
\begin{eqnarray*}
{{\mathcal  L}}_{n} u_n &=& g,
\\
u_n (T,\cdot) &=& u^0.
\end{eqnarray*}
According to points (ii) and (vii) of the same lemma, one can prove
that $u_n$ converges in a bounded way to $u$ defined by
(\ref{F312terd}). In fact, the coefficients $\sigma^2_n$ are lower and
upper bounded with a common constant, related to $c$ and $C$.
Therefore, this $u$ is the $C^0_b$-generalized solution of the Cauchy
problem being considered, which is known to exist. By uniqueness,
point~(i) is established.

\item[(ii)] Setting $g = 0$, point (i) implies that $ u (s,x) =
\int_{\mathbb R} \varphi(y) r_{T-s} (x,y)\,dy $ is the
$C^0_b$-generalized solution to ${\mathcal  L}u = 0$ with $ u(T,x) =
\varphi(x) $. By Proposition~\ref{Probrep}, in particular, using the
probabilistic representation, we get\break ${\mathbb E} (\varphi(Z_T^{s,x}) )
= \int_{\mathbb R} \varphi(y) r_{T-s} (x,y)\,dy $.\quad\qed
\end{longlist}\noqed
\end{pf}

\begin{rem} \label{Rdistr}
If $L$ is in the divergence form, as before, then ${\mathcal  D}_{L} =
\{ f \in C^1$ such that there exists $g \in C^1 $ with $ f' = \frac{g}
{\sigma^2} \} $. This is a consequence of Lemma \ref{L28bis} and the
fact that $e^{-\Sigma} = \frac{1}{\sigma^2} $.
\end{rem}

{\renewcommand{\theequation}{\textup{Aronson}}
Hereafter, we will
consider a general PDE operator $L$ with distributional drift, as in
Section~\ref{s4}, for which the assumption~(\ref{Aronson}) below holds.
\begin{equation}\label{Aronson}
c \le\frac{e^\Sigma} {\sigma^2} \le C. 
\end{equation}}%
\noindent We observe that the PDE operator in divergence form of the type $L^1 f
= (\frac{\sigma_k ^2 f'}{2})'$, where $\sigma_k = (\sigma k')
\circ{k^{-1}}$, has the Aronson form, so the previous theorem can be
applied.

\setcounter{equation}{3}
\begin{theo} \label{TDistr1}
Let $X^{s,x}$ be the solution to the martingale problem related to $L$
at time $s$ and point $x$. Suppose that $L$ fulfills
assumption~\textup{(\ref{Aronson})}. Then there exists a kernel
$p_t(x,y)$ such that:
\begin{longlist}[(ii)]
\item[(i)] the law of $X_t^{s,x}$ has $p_{t-s}(x,\cdot)$ as density
with respect to Lebesgue measure for each $t \in\,]s,T]$;

\item[(ii)] letting $g \in C^0_b([0,T] \times{\mathbb R})$, $\varphi\in
C^0_b({\mathbb R})$, the $C_b^0$-generalized solution $u$ to ${\mathcal
L}u=g$, $ u(T,\cdot) =\varphi$, is given by
\begin{equation} \label{F312tera}
\qquad u(s,x) = \int_{\mathbb R}\varphi(y) p_{T-s}(x,y)\,dy + \int_s^T\,dr
\int_{\mathbb R}g(r,y) p_{T-r} (x,y)\,dy .
\end{equation}
\end{longlist}
\end{theo}

\begin{pf}
(i) Proposition~\ref{PMardiv} says that $Z^{s,x} = k(X^{s,x})$ solves
the martingale problem with respect to $L^1$. Let $r_t(x,y)$ be the
fundamental solution associated with the parabolic PDE ${\mathcal  L}^1
= \partial_t + L^1$. The first point then follows from the next
observation.

\begin{rem}\label{RCV}
By means of a change of variable, it is easy to see that the density
law of $X_t^{s,x}$ equals
\[
p_t(x,x_1) = r_t(k(x),k(x_1)) k'(x_1) = r_t(k(x),k(x_1)) \frac{
e^{\Sigma}}{\sigma^2} (x_1).
\]
\end{rem}

(ii) This is a consequence of point (i), Fubini's theorem and
Proposition~\ref{Probrep}.
\end{pf}

At this point, we need a lemma which extends to the kernel $p_t(x,x_1)$
the integrability property of the kernel $r_t(x,x_1)$ stated in
(\ref{Rdistr}) concerning the divergence case.

\begin{lemma} \label{DDistr1}
Let $p_t(x,x_1)$ be the kernel introduced in
Theorem~\textup{\ref{TDistr1}}. Then:
\begin{longlist}[(iii)]
\item[(i)] it is continuous in all variables $(t,x,x_1) \in\,]0,T[\,
\times\,
{\mathbb R}^2$;

\item[(ii)] it fulfills Aronson estimates;

\item[(iii)] $ \int_0^T (\sup_{x_1} \int_{\mathbb R} \partial_x
p_t(x,x_1)^2\,dx )^{{1}/{2}}\,dt < \infty$.
\end{longlist}
\end{lemma}

\begin{pf}
We recall, by Remark~\ref{RCV}, that
\[
p_t(x,x_1) = r_t(k(x),k(x_1)) k'(x_1),
\]
where $r_t(z,z_1)$ is the fundamental solution associated with the
operator $L^1 f = (\frac{ \sigma_k^2}{2} f')'$,  $k' = \frac{e^{\Sigma
}}{\sigma^2} $. This, and point (v) of Lemma \ref{L5.1}, directly imply
the validity of the first point.

Taking into account assumption~\textup{(\ref{Aronson})}, Aronson
estimates for $(r_t(z,z_1))$ and the fact that
\[
| k(x) -k(x_1)|= \int_0^1 k'\bigl(\alpha x + (1-\alpha)
x_1\bigr)\,d\alpha| x - x_1 |,
\]
result (ii) follows easily.

With the same conventions as before, we have
\[
\partial_x p_t(x,x_1) = \partial_z r_t(k(x),k(x_1)) k'(x) k'(x_1).
\]

So, for $x \in{\mathbb R}$,
\begin{eqnarray*}
 \biggl( \int_{\mathbb R} (\partial_x p_t (x,x_1))^2 \,dx  \biggr)
^{{1}/{2}} &=& \biggl( k'(x_1) \int_{\mathbb R} (\partial_z r_t
(z,k(x_1)))^2 \,dz \biggr) ^{{1}/{2}}
\\
&\le& \sqrt C \sup_{z_1}  \biggl( \int_{\mathbb R}dz (\partial_z
r_t(z,z_1))^2  \biggr)^{{1}/{2}}.
\end{eqnarray*}

(iii) Follows after integration with respect to $t$ and because of
Lemma~\ref{L5.1}(ix).
\end{pf}

\begin{prop}\label{PDistr1}
Let $g \in C^0_b([0,T] \times{\mathbb R}) \cap L^1([0,T] \times
{\mathbb R})$, $ \varphi\in C^0_b({\mathbb R}) \cap\break L^1(
{\mathbb R}) $. Let $u \dvtx  [0,T] \times{\mathbb
R}\rightarrow{\mathbb R}$ be the $C^0_b$-generalized solution to
${\mathcal  L}u = g$,\break $u(T, \cdot) = \varphi$. Then:
\begin{longlist}[(a)]
\item[(a)] $ \int_0^T dt \int_{\mathbb R}u^2(t,x)\,dx < \infty$;

\item[(b)] $x \mapsto u(t,x) $ is absolutely continuous,
\[
\int_0^T dt  \biggl(\int_{\mathbb R}(\partial_x u)^2(t,x)\,dx
\biggr)^{1/2}< \infty
\]
and in particular, for a.e. $t \in[0,T]$, $\partial_x u(t, \cdot)$ is
square integrable.
\end{longlist}
\end{prop}

\begin{rem} \label{RDistr1}
Previous assumptions imply that $g$ and $\varphi$ are also square integrable.
\end{rem}

\begin{pf*}{Proof of Proposition~\ref{PDistr1}}
We recall the expression given in Theorem \ref{TDistr1},
\[
u(t,x) = \int_{\mathbb R}\varphi(x_1) p_{T-t}(x,x_1)\,dx_1 +\int_t^T
dr \int_{\mathbb R}g(r,x_1) p_{T-r} (x,x_1)\,dx_1.
\]
Using Lemma \ref{DDistr1} and classical integration theorems, we have
\begin{eqnarray} \label{deru}
\partial_x u (t,x) &=& \int_{\mathbb R}\varphi(x_1)\, \partial_x
p_{T-t}(x,x_1)\,dx_1
\nonumber
\\[-8pt]
\\[-8pt]
\nonumber &&{} + \int_t^T dr \int_{\mathbb R} \,ds\, g(s,x_1)\,
\partial_x p_{T-s} (x,x_1)\,dx_1 .
\end{eqnarray}
Using Jensen's inequality, we have
\begin{eqnarray*}
| u(t,x) |^2 &\le& \int_{\mathbb R} \varphi(x_1)^2 p_{T-t}(x,x_1)\,dx_1
\\
&&{} + (T-t) \int_t^T ds \int_{\mathbb R}g^2 (s,x_1) p_{T-s}
(x,x_1)\,dx_1 .
\end{eqnarray*}
Therefore,
\begin{eqnarray*}
\int_{\mathbb R}u^2 (t,x)\,dx &=& \int_{\mathbb R}dx_1\, \varphi(x_1)^2
\int_{\mathbb R}dx\, p_{T-t}(x,x_1)
\\
&&{} + \int_t^T ds\, (T-t) \int_{\mathbb R} dx_1 \int g^2 (s,x_1) \int
_{\mathbb R}dx\, p_{T-s} (x,x_1).
\end{eqnarray*}
Using Aronson estimates, this quantity is bounded by
\begin{eqnarray*}
&& \operatorname{const}  \biggl(\int_{\mathbb R} dx_1\, \varphi(x_1)^2
\int _{\mathbb R}dx\, \frac{1}{\sqrt{T-t}} p
\biggl(\frac{x-x_1}{\sqrt{T-t}} \biggr)
\\
&& \phantom{\operatorname{const}  \biggl(} {}  + \int_t^T ds
\int_{\mathbb R}dx_1 \int g^2 (s,x_1) \int _{\mathbb R}dx
\,\frac{1}{\sqrt{T-s}} p \biggl(\frac{x-x_1}{\sqrt{T-s}}
\biggr)\biggr),
\end{eqnarray*}
where $p$ is the Gaussian $N(0,1)$ density.
This is clearly equal to
\[
\operatorname{const}  \biggl(\int_{\mathbb R}dx_1\, \varphi(x_1)^2 +
\int_0^T ds \int_{\mathbb R}dx_1\, g^2 (s,x_1)  \biggr).
\]
This establishes point \textup{(a)}.

Concerning point (b), in order to simplify the framework we will
suppose that $g = 0$. Expression (\ref{F312tera}) implies that
\[
\partial_x u(t,x) =
\int_{\mathbb R} \varphi(x_1)\, \partial_x p_{T-t} (x,x_1)\,dx_1 .
\]
Jensen's inequality implies that
\[
\partial_x u(t,x)^2 \le \biggl( \int_{\mathbb R} dx_1\, |\varphi
(x_1)||
\partial_x p_{T-t} (x,x_1)|^2\biggr)   \int_{\mathbb R} dx_1\,
|\varphi(x_1)|.
\]
Integrating with respect to $ x$ and taking the square root, we get
\begin{eqnarray*}
\sqrt{ \int_{\mathbb R}dx \,\partial_x u(t,x)^2} &\le& \sqrt{ \int
_{\mathbb R}|\varphi(x_1) |\,dx_1 } \sqrt{ \int_{\mathbb R} dx_1\,
|\varphi(x_1) | \int_{\mathbb R}dx \,|\partial_x p_{T-t} (x,x_1)|^2}
\\
&\le& \int_{\mathbb R} dx_1\, |\varphi(x_1) | \sqrt{ \sup_{x_1}
\int_{\mathbb R} |\partial_x p_t(x,x_1) |^2\, dx}.
\end{eqnarray*}
Integrating with respect to $t$ gives
\[
\int_0^T dt\, \Vert\partial_x u(t,\cdot) \Vert_{L^2({\mathbb R})}
\le\int_{\mathbb R} dx_1 \,|\varphi(x_1) | \int_0^T dt \sqrt{\sup_{x_1}
\int_{\mathbb R} \partial_x p_t(x,x_1)^2\, dx}.
\]
This quantity is finite due to Lemma \textup{\ref{DDistr1}(iii)}.
\end{pf*}

\section{Relation with weak solutions of stochastic partial
differential equations}\label{s9}

As in the previous section, we will adopt
assumption~\textup{(\ref{Aronson})}. At this point, we wish to
investigate the link between $C^0_b$-generalized solutions and the
notion of SPDE's weak solutions for a corresponding Cauchy problem.

We will adopt the same conventions as in Section~\ref{s4}. In this
section, we will suppose that coefficients $\sigma, b$ are realizations
of stochastic processes indexed by~${\mathbb R}$. Let us consider the
formal operator ${\mathcal  L}= \partial_t + L$, where $L$ acts on the
second variable.

We consider the equation
\begin{eqnarray} \label{E3.9quater}
{\mathcal  L}u &=& \lambda,
\nonumber
\\[-8pt]
\\[-8pt]
\nonumber u(T,\cdot) &=& u^0.
\end{eqnarray}

The aim of this section is to show that a $C^0_b$-generalized solution
to (\ref{E3.9quater}) provides, when $\sigma= 1$, a solution to the
(stochastic) PDE of the type (\ref{eI-0}), as defined in the
Definition \ref{deI-1rig}, that is, with the help of a symmetric
integral via regularization, as defined in Section~\ref{s3}. We denote
by ${\mathcal D}({\mathbb R})$ the linear space of $C^\infty$ real
functions with compact support.

The link between the SPDE (\ref{eI-0}) and (\ref{eI-0dual})
is given in the following.

\begin{prop} \label{PSPDEs}
Let $u(t,x)$, $v(t,x)$, $t \in[0,T]$, $x \in{\mathbb R}$ be two
continuous random fields a.s. in $C^{0,1} (]0,T[ \times{\mathbb R}) $
such that $v(t,x) = u(T-t,x)$. $v$ is a solution to the SPDE
\textup{(\ref{eI-0})} if and only if $v$ is a solution to the SPDE
\textup{(\ref{eI-0dual})}.
\end{prop}

\begin{pf}
We observe that $ \partial_x v(t,x) = - \partial_x u(T-t,x)$.
The proof is elementary. The only point to check is the following:
\[
\int_{\mathbb R}d^\circ\eta(x) \alpha(x)  \biggl( \int_t^T ds\,
\partial_x u(s,x)  \biggr)
= - \int_{\mathbb R}d^\circ\eta(x) \alpha(x)  \biggl( \int_0^t ds\,
\partial_x v(s,x)  \biggr).
\]
This follows by the definition of symmetric integral and the following,
obvious, identity:
\begin{eqnarray*}
&& \int_{\mathbb R} dx\, \frac{\eta(x + \varepsilon) - \eta(x -
\varepsilon)}{2 \varepsilon} \alpha(x)  \biggl( \int_t^T ds
\,\partial_x u(s,x)
 \biggr)
\\
&&\qquad = - \int_{\mathbb R} dx\, \frac{\eta(x + \varepsilon) - \eta(x -
\varepsilon)}{2 \varepsilon} \alpha(x)  \biggl( \int_0^t ds\, \partial_x
v(s,x) \biggr)
\end{eqnarray*}
for every $\varepsilon> 0$.
\end{pf}

We continue with a lemma, still supposing $\sigma$ to be general.

\begin{lemma} \label{Lspdes}
Let $\lambda$ (resp., $u^0$) be a random field with parameter $(t,x)
\in[0,T] \times{\mathbb R}$ (resp., $x \in{\mathbb R}$) whose paths are
bounded and continuous. Let $\sigma, b$ be continuous stochastic
processes such that $\Sigma$ is defined a.s. and
assumption~\mbox{\textup{(\ref{Aronson})}} is satisfied. Let $u$ be the
random field which is a.s. the $C^0_b$-generalized solution
to~(\ref{E3.9quater}). The following then holds:
\begin{eqnarray*}
&& \int_{\mathbb R}dx\, \alpha(x)  \biggl(u(t,x) - u^0(x) + \int_t^T
\lambda(s,x)\,ds \biggr)
\\
&&\qquad = \int_{\mathbb R} e^{\Sigma(x)}  \biggl( \int_t^T ds\,
\partial_x u(s,x)
 \biggr)\,d^\circ \biggl(\alpha\frac{\sigma^2}{2} e^{-\Sigma}(x)  \biggr)
\end{eqnarray*}
\end{lemma}
for every $\alpha\in{\mathcal  D}({\mathbb R})$.

\begin{pf}
We fix a realization $\omega$. Theorem \ref{TDistr1} says that the
unique solution to equation (\ref{E3.9quater}) is given by
\begin{equation} \label{F312terb}
\qquad u(s,x) = \int_{\mathbb R}u^0 (y) p_{T-s}(x,y)\,dy + \int_s^T dr\,
\int_{\mathbb R}\lambda(r,y) p_{T-r} (x,y)\,dy,
\end{equation}
where $(p_t(x,y))$ is the density law of the solution to the martingale
problem related to $L$ at point $x$ at time $s$.

Proposition~\ref{PDistr1}(b) implies that $\partial_x u$ exists and
is integrable on $]0,T[ \times{\mathbb R}$.

According to Proposition~\ref{P313}, we know that
\[
\bar u(t, z) = u(t, k^{-1}(z))
\]
is a $C^0_b$-generalized solution to
\begin{eqnarray}\label{C0div}
{{\mathcal  L}}^1 \bar{u} &=& \bar\lambda,
\nonumber
\\[-8pt]
\\[-8pt]
\nonumber \bar{u} (T,\cdot) &=& \bar u^0,
\end{eqnarray}
where
\[
\bar\lambda(t, z) = \lambda(t, k^{-1}(z)),\qquad \bar u^0 (z) = u^0
(k^{-1}(z)).
\]
On the other hand, $\bar u$ can be represented via (\ref{F312terd}) in
Theorem \ref{TDistr} through fundamental solutions $(\nu_t)
=(r_t(x,y))$ of
\[
\partial_t \nu_t (\cdot,y) = L^1 \nu_t(\cdot,y),\qquad
\nu_0(\cdot,y) = \delta_y.
\]
Since the previous equation holds in the Schwarz distribution sense,
by inspection, it is not difficult to show that $\bar u$ is a solution
(in the sense of distributions) to (\ref{C0div}), which means that we
have the following:
\begin{eqnarray} \label{C0div1}
&& \int_{\mathbb R}\alpha(z) \bigl(\bar u^0 (z) - \bar u(t,z)\bigr)\,dz
- \tfrac {1}{2}\int_t^T ds \int_{\mathbb R}\alpha'(z)\, \partial_z \bar
u(s,z) \sigma_k^2 (z) \nonumber
\\[-8pt]
\\[-8pt]
\nonumber &&\qquad = \int_t^T ds \int_{\mathbb R}\alpha(z)
\bar\lambda(s,z)
\end{eqnarray}
for every test function $\alpha\in{\mathcal  D}({\mathbb R})$, $t
\in[0,T]$. We recall, in particular, that $\partial_z \bar u$ is in
$L^1(]0,T[ \times{\mathbb R})$.

We set
\[
D(t,z) = \int_t^T \partial_z \bar u(s,z)\,ds,\qquad  {\mathcal  D}(t,z)
= D(t,z) \frac{\sigma_k^2(z)}{2}.
\]
Expression (\ref{C0div1}) shows that
\begin{equation} \label{EDiv1}
\partial_z {\mathcal  D}(t,\cdot) = - \bar u^0 + \bar u(t, \cdot) + \int
_t^T \bar\lambda(s,\cdot)\,ds,
\end{equation}
in the sense of distributions. So for each $t \in[0,T]$, ${\mathcal D}$
is of class $C^1$.

For $t \in[0,T]$ and $x \in{\mathbb R}$, we set $A(t,x) = \int_t^T
\partial_x u(s,x)\,ds$, ${\mathcal  A}(t,x) =\break A(t,x) e^{\Sigma(x)}$. We
recall that
\[
u(s,x) = \bar u(s,k(x)),\qquad  \partial_x u(s,x) = \partial_x \bar
u(s,k(x)) k'(x).
\]
Therefore,
\[
A(t,x) = D(t,k(x)) k'(x)
\]
so that
\begin{eqnarray*}
A(t,x) &= & 2 {\mathcal  D} (t,k(x)) \frac{k'(x)}{\sigma_k^2(k(x))} =
{\mathcal  D}(t,k(x))\frac{2}{\sigma^2(x) k'(x) }
\\
&=& 2 {\mathcal  D} (t,k(x)) e^{-\Sigma(x)}.
\end{eqnarray*}
Therefore, $ {\mathcal  A}(t,x) = 2 {\mathcal  D}(t,k(x)) $ and so
${\mathcal  A}$ is of class $C^1$.

Since $ \partial_x {\mathcal  A}(t,x) = 2 \partial_z {\mathcal
D}(t,k(x)) k'(x) $, (\ref{EDiv1}) gives
\begin{equation} \label{EPoint3}
\partial_x {\mathcal  A}(t,x) =  \biggl( - u^0(x) + u(t, x) + \int_t^T
\lambda(s,x)\,ds  \biggr) 2 \frac{e^{\Sigma}}{\sigma^2}(x).
\end{equation}

Consequently,
\[
u(t,x) - u^0 (x) + \int_t^T \lambda(s,x)\,ds =
\partial_x  \biggl( e^{ \Sigma(x)} \int_t^T ds \,\partial_x u(s,x)
\biggr) e^{- \Sigma(x)}\frac{\sigma^2 (x)}{2}.
\]
We integrate the previous expression against a test function $\alpha\in
{\mathcal  D}({\mathbb R})$ to obtain
\begin{eqnarray*}
&& \int_{\mathbb R}dx\, \alpha(x)  \biggl(u(t,x) - u^0 (x) + \int_t^T
\lambda(s,x)\,ds  \biggr)
\\
&&\qquad = \int_{\mathbb R}dx\, \alpha(x)  \biggl\{ \partial_x  \biggl(
e^{ \Sigma (x)} \int_t^T ds\, \partial_x u(s,x)  \biggr) e^{- \Sigma(x)}
\frac{\sigma^2 (x)}{2} \biggr\}.
\end{eqnarray*}

Remark~\ref{R1.0} and integration by parts for the symmetric integral
provided by Remark~\ref{R1.1}(c) allow us to conclude the proof of the
lemma.
\end{pf}

Finally, we are able to state the theorem concerning the existence of
weak solutions for the SPDE.

\begin{theo} \label{Tspdes}
Let $\lambda$ (resp., $u^0$) be a random field with parameter in $(t,x)
\in[0,T] \times{\mathbb R}$ (resp., $x \in{\mathbb R}$) whose paths are
bounded and continuous.

We suppose that $\sigma=1$ and that $\eta$ is a (two-sided) zero strong
cubic variation process such that there are two finite and strictly
positive random variables $Z_1, Z_2$ with $ Z_1 \le e^{\eta(x)}
\le Z_2 $ a.s.

Let $u$ be the random field which is $\omega$ a.s. a
$C^0_b$-generalized solution to~\textup{(\ref{E3.9quater})} for $b =
\eta(\omega)$. We set $v(t,x) = u(T-t,x)$. Then $v$ is a (weak)
solution of the SPDE~\textup{(\ref{eI-0})}.
\end{theo}

\begin{pf}
Proposition~\ref{PSPDEs} says that it will be enough to verify that
\begin{eqnarray*}
&& - \int_{\mathbb R}\alpha(x) u(t,x)\,dx + \int_{\mathbb R}\alpha(x)
u^0(x)\,dx
\\
&&\quad {} - \tfrac{1}{2}\int_{\mathbb R}\alpha'(x)  \biggl(\int_t^T ds
\,\partial _x u(s,x)  \biggr)\,dx +  \int_{\mathbb R}\alpha(x)
\biggl(\int_t^T ds \,\partial_x u(s,x) \biggr)\,d^\circ\eta(x)
\\
&&\qquad = \int_t^T ds \int_{\mathbb R}dx\, \alpha(x) \lambda(s,x)
\end{eqnarray*}
for every test function $\alpha$ and every $t \in[0,T]$.

After making the identification $b = \eta(\omega)$, the previous Lemma
\ref{Lspdes} says that
\begin{eqnarray*}
&& \int_{\mathbb R} dx\, \alpha(x)  \biggl(u(t,x) - u^0(x) + \int_t^T
\lambda(s,x)\,ds  \biggr)
\\
&&\qquad =  \int_{\mathbb R} e^{2\eta(x)}  \biggl( \int_t^T ds\,
\partial_x u(s,x)\biggr) \,d^\circ \biggl( \frac{\alpha e^{-2 \eta}}{2} (x)
\biggr).
\end{eqnarray*}

Since $\eta$ is a zero strong cubic variation process,
Proposition~\ref{Pstacub} implies that $e^\eta$ is also a zero strong
cubic variation process. Then the It\^o chain rule from
Proposition~\ref{Crule}, applied with $ F(x, \eta(x)) = \alpha(x)
e^{\eta(x)},$ and Remark~\ref{R1.0} say that the right-hand side of
previous expression gives
\begin{eqnarray*}
&& \tfrac{1}{2}\int_{\mathbb R} \biggl( \int_t^T ds \,\partial_x u(s,x)
\biggr)\,d^{0} \bigl(\alpha e^{-2 \eta(x)}\bigr)
\\
&&\qquad = \tfrac{1}{2}\int_{\mathbb R}
 \biggl(\int_t^T ds \,\partial_x u(s,x)  \biggr) e^{2 \eta(x)}
\bigl(\alpha'(x) e^{-2 \eta}(x)\,dx + \alpha(x)\,d^\circ e^{- 2 \eta(x)
}  \bigr)
\\
&&\qquad = \tfrac{1}{2}\int_{\mathbb R}  \biggl(\int_t^T ds\,
\partial_x u(s,x)  \biggr) \alpha'(x)\,dx
\\
&&\qquad\quad{} - \int_{\mathbb R}  \biggl(\int_t^T ds\, \partial_x
u(s,x) \biggr) \alpha(x)\,d^\circ\eta(x).
\end{eqnarray*}
This concludes the proof.
\end{pf}

\section*{Acknowledgments}
We  would like to thank an anonymous referee  and the Editor
for their careful reading and stimulating remarks. The first named
author is grateful to Dr. Juliet Ryan for her precious help in
correcting several language mistakes.

\printaddresses

\end{document}